\documentclass[11pt]{article}
\usepackage{mathrsfs}
\usepackage{amssymb}
\usepackage{amsmath}
\usepackage{amsbsy}
\usepackage{epsfig}
\usepackage{enumerate}
\usepackage{bm}
\newtheorem{lemma}{Lemma}[section]
\newtheorem{theorem}{Theorem}[section]
\newtheorem{proposition}{Proposition}[section]
\newtheorem{assumption}{Assumption}[section]

\newtheorem{remark}{Remark}[section]

\newbox\TempBox \newbox\TempBoxA

\def\ep{\textsf{E}} 
\def\Cov{\textsf{Cov}} 
\def\Var{\textsf{Var}} 
\def\Cal#1{{\mathcal #1}}
\def\bk#1{\mathbf #1}

\def\underwiggle 1{
\ifmmode\setbox\TempBox=\hbox{$ 1$}\else\setbox\TempBox=\hbox{
1}\fi \setbox\TempBoxA=\hbox to \wd\TempBox{\hss\char'176\hss}
\rlap{\copy\TempBox}\smash{\lower9pt\hbox{\copy\TempBoxA}} }
\renewcommand{\baselinestretch}{1.7}

\begin{document}

\thispagestyle{empty}

\begin{center}
{ \LARGE\bf The Gaussian approximation for multi-color generalized
Friedman's urn model $^{\ast}$}
\end{center}

\begin{center} {\sc
Li-Xin Zhang\footnote{Research supported by grant 10771192 from the
National Natural Science Foundation of China.},
 Feifang
Hu\footnote{Research supported by grant   DMS-0349048 from the
National Science Foundation (USA).}}

{\sl \small Zhejiang University and University of Virginia}
\end{center}

\renewcommand{\abstractname}{~}
\begin{abstract}
\centerline{\bf Abstract}
 The   Friedman's urn model is a popular urn model which is widely used in many disciplines. In particular, it is extensively
 used in treatment allocation schemes in clinical trials.
 In this paper, we prove that both the urn composition process and the allocation proportion process
   can be approximated by a multi-dimensional Gaussian
 process  almost surely
 for a   multi-color generalized
Friedman's urn model  with non-homogeneous generating matrices.
 The  Gaussian process is a solution of a stochastic differential
 equation. This Gaussian approximation  together with the properties
 of the Gaussian process is important for the understanding of the behavior of the
 urn process and is also useful  for statistical inferences. As an
 application, we obtain the asymptotic properties including the
 asymptotic normality and the law of the iterated logarithm for a multi-color generalized
Friedman's urn model as well as the randomized-play-the-winner rule
as a special case.

 \bigskip
{\bf Abbreviated Title:} Approximation for multi-color urn models
\bigskip

\bigskip


\end{abstract}

\baselineskip 22pt
\newpage
\renewcommand{\baselinestretch}{1.7}



\section{ Introduction.}
\setcounter{equation}{0}

Urn models have long been recognized as valuable mathematical
apparatus in many areas including  physical sciences, biological
sciences and engineering (Johnson and Kotz, 1977; Kotz and
Balakrishnan, 1997). Urn models are also extensively applied in
clinical studies.   The applications are mostly found in the area of
adaptive design which is utilized to provide a response-adaptive
allocation scheme.  In  clinical trials, suppose patients accrue
sequentially and  assume the availability of several treatments.
Adaptive designs are inclining to assign more patients to the better
treatments, while seeking to maintain randomness as a basis for
statistical inference. Thus the cumulative information of the
response of treatments on previous patients will be used to adjust
treatment assignment to coming patients. For this purpose, various
urn models    have been proposed and used
extensively in adaptive designs (Wei and Durham (1978), Wei (1979),
Flournoy and Rosenberger (1995), Rosenberger (1996), Bai and Hu
(1999,2005)).
 One large family of randomized adaptive designs is
based on the Generalized Friedman's Urn (GFU) model (also named as
Generalized P\'olya Urn (GPU) in literature). For more detailed
reference, the reader is referred to Flournoy and Rosenberger
(1995), Rosenberger (1996), Hu and Rosenberger (2006).

A general description of the GPU model is as follows.   Consider an
urn containing particles of $d$ types, respectively representing $d$
'treatments' in a clinical trial. At the beginning, the urn contains
$\bm Y_0=(Y_{01}, \ldots, Y_{0d})$ particles, where $Y_{0k}>0$
denotes the number of particles of type $k$, $k=1,\ldots, d$. At the
stage $m$, $m=1,2,\ldots$, a particle is drawn form the urn and
replaced. If the particle is of type $k$, then the treatment $k$ is
assigned to the $m$th patient, $k=1,\ldots, d$. We then wait for
observing a random variable $\bm \xi(m)$, the response of the
treatment at the patient $m$. Here $\bm\xi(m)$ may be a random
vector. After that, an additional $D_{k,q}(m)$ particles of type
$q$, $q=1,\ldots, d$, are added to the urn, where $D_{k,q}(m)$ is a
function of $\bm \xi(m)$ and also may be a function of urn
compositions, assignments and responses of previous stages. This
procedure is repeated through out $n$ stages. After $n$ draws and
generations, the urn composition is denoted by the row vector $\bm
Y_n=(Y_{n1}, \ldots, Y_{nd})$, where $Y_{nk}$ stands for the number
of particles of type $k$ in the urn after the $n$th draw. This
relation can be written as the following recursive formula:
\begin{equation}\label{eq1.1}
 \bm Y_m = \bm Y_{m-1}+\bm X_m \bm D_m, \end{equation}
 where $\bm D_m=\big(D_{k,q}(m)\big)_{k,q=1}^d $, and $\bm X_m$ is the result
of the $n$th draw, distributed according to the urn composition at
the previous stage, i.e., if the $m$th draw is type $k$ particle,
then the $k$th component of $\bm X_m$ is $1$ and other components
are $0$. The matrices $\bm D_m$'s are named as the {\bf addition
rules}. Furthermore, write $\bm N_n=(N_{n1},\ldots, N_{nd})$, where
$N_{nk}$ is the number of times a type $k$ particle drawn in the
first $n$ stages. In clinical trials, $N_{nk}$ represents the number
of patients assigned to the treatment $k$ in the first $n$ trials.
Obviously,
\begin{equation}\label{eq1.2}
 \bm N_n=\sum_{k=1}^n \bm X_n. \end{equation}
 In clinical applications,  $\bm Y_{n-1}/\sum_{k=1}^d Y_{n-1,k}$ are the
 probabilities of the patient $n$ being allocated to treatments,
 and $\bm N_n/n$ are sample allocation proportions.
The asymptotic behavior of $\bm Y_n$ and  $\bm N_n$ is of immense
importance (Hu and Rosenberger, 2003, 2006). Obviously, the
asymptotic behavior of $\bm Y_n$ and  $\bm N_n$ will depend on the
addition  rules $\bm D_m$, especially the conditional expectations
$\bm H_m=\big(\ep[D_{k,q}(m)\big|\mathscr{F}_{m-1}]\big)_{k,q=1}^d $
for given the history sigma field $\mathscr{F}_{m-1}$ generated by
the urn compositions  $\bm Y_1,\ldots, \bm Y_{m-1}$, the assignments
$\bm X_1,\ldots, \bm X_{m-1}$ and the responses $\bm \xi(1),\ldots,
\bm \xi(m-1)$ of all previous stages, $m=1,2,\ldots$. The
conditional expectations $\bm H_m$'s are named as the {\bf
generating matrices}. In some usual cases, the addition rules are
assumed to be independent of the previous process. Thus, we may
define $\bm H_m$ is the expectation of the rule matrix $\bm D_m$.
For more generality, in the sequel of this paper, we define $\bm
H_m$ to be the conditional expectation of $\bm D_m$ when the history
sigma field $\mathscr{F}_{m-1}$ is given, also we assume that at the stage $m$,  the adding rule $\bm
D_m$ is independent of the assignment $\bm X_m$ when the history sigma field is
given.

 When $\bm
D_m$, $m=1,2,\ldots, $ are independent and identical distributed,
the GFU model is usually  said to be {\bf homogeneous}. In such case
$\bm H_m=\bm H$ are identical and nonrandom, and usually the
addition rule $\bm D_m$ is merely function of the $m$th patient's
observed outcome. In the general non-homogeneous cases, both $\bm
D_m$ and $\bm H_m$ depend on the entire history of all previous
trials which provides more information of the efficacy of the
treatments. Interesting examples of non-homogeneous urn models and
their applications can be found in Andersen, Faries and Tamura
(1994) and Bai, Hu and Shen (2002).

 Athreya and
Karlin (1967, 1968) first considered the asymptotic properties of
the GFU model with homogeneous generating matrix and conjecture that
$\bk N_n$ is asymptotically normal.  This conjecture has not been
solved for almost three decades  until  Janson (2004) and Bai and Hu
(2005) solved it independently. Janson (2004) established functional
limit theorems of $\bk Y_n$ and $\bk N_n$ for the  homogenous case
by using the theory of continuous-time branching processes. Bai and
Hu (2005) established the consistency and the asymptotic normality
of the  non-homogeneous GFU model by applying the central limit of
martingales and the matrix theory. However, the  asymptotic
variances of $\bk Y_n$ and $\bk N_n$ are complicated  and not easy
to be understood.

In the two-arm clinical trial,  Bai, Hu and Zhang (2002) showed that
the urn process $\{Y_n\}$ with nonhomogeneous generating matrices
$\bm H_m$'s can be approximated by a Gaussian process almost surely
under some suitable conditions,  where $Y_n=Y_{n1}$ represents the
number of type 1 balls in the urn after the $n$th draw. As an
application, the weak invariance principle and the law of the
iterated logarithm for $\{Y_n\}$ are established. However, the
results for the allocation proportion $N_{n1}/n$  is not obtained.
In this paper, we consider the general multi-color case. The strong
approximation of the process $(\bm Y_n, \bm N_n)$ are established.
In particular, the asymptotic normality and the law of the iterated
logarithm for the multi-dimensional process $(\bm Y_n, \bm N_n)$ are
obtained. We will prove that under some mild conditions, the process
$(\bm Y_n, \bm N_n)$ can be approximated by a multi-dimensional
Gaussian process which is a solution of a simple multi-dimensional
stochastic differential equation. This  differential equation and
the behavior of the Gaussian process make us  to understand the
complex asymptotic variances and the asymptotic behavior of $\bm
Y_n$ and $\bm N_n$ more easily.

The approximation theorems will be presented in Section 2 whose
technical proofs are stated in the last section. Some important
properties of the limit processes are given in Section 3. By
combining the approximation theorems and  the properties of the
limit processes, important asymptotic properties including the
asymptotic normality and the law of the iterated logarithm  of $\bm
Y_n$ and $\bm N_n$ are derived in Section 4. Throughout this paper,
$C, C_{\epsilon}, etc.$ denote positive constants whose values can
differ in different places, $\log x=\ln(e\vee x)$. For a vector $\bm
x$, $\|\bm x\|$ denote its Euclidean norm, and  $\|\bm
M\|=\sup\{\|\bm x\bm M\|: \|\bm x\|=1\}$ for a martix $\bm M$. Also
we denote $a_n=\bm Y_n \bm 1^{\prime}$ to be the total number of
balls in the urn after stage $n$.

\section{Strong approximation.}
\setcounter{equation}{0} In this section, we will give our main
results on the Gaussian approximation of both the urn composition
$\bm Y_n$ and the allocation numbers $\bm N_n$. We first need two
assumptions on the addition  rules $\bm D_m$. We let $\mathscr{F}_m=
\sigma(\bm Y_1,\ldots, \bm Y_m, \bm X_1,\ldots,\bm X_m,\bm
\xi(1),\ldots,\bm \xi(m))$ be  history sigma field. and $\bm
H_m=\ep[\bm D_m|\mathscr{F}_{m-1}]$ be the generating matrix.

\begin{assumption} \label{assumption1.1} Suppose there is a $\tau\ge 0$ such that the generating matrices $\bm H_m$
satisfy
\begin{equation}\label{eqassumption1.1.1}
\sum_{m=1}^n \|\bm H_m-\bm H\|=o(n^{1/2-\tau}) \;\; a.s.,
\end{equation}
with $\bm H$ bing  a deterministic matrix and
$$H_{qk}\ge 0 \;\text{ for }\; k\ne q \; \text{ and } \; \sum_{k=1}^d H_{qk}=s
\; \text{ for all}\; q=1,\ldots, d, $$ where $H_{qk}$ is the
$(q,k)$-entry of the matrix $\bm H$ and $s$ is a positive constant.
Without loss of generality, we assume $s=1$ through out this paper.
For otherwise, we may consider $\bm Y_m/s$, $\bm H_m/s$ instead.
\end{assumption}

\begin{assumption} \label{assumption1.3} Let
$$ V_{qkl}(n)=:\Cov\big[(D_{qk}(n), D_{ql}(n))\big|\mathscr{F}_{n-1}\big], \quad q,k,l=1,\ldots d $$
and denote by $\bm V_{nq}= (V_{qkl}(n))_{k,l=1}^d$.
Suppose for
some $0<\epsilon<1/2$,
\begin{equation}\label{eq1.4}\ep \|\bm D_n\|^{2+\epsilon}\le C<\infty \; \text{ for all } \; n=1,2,
\ldots,  \end{equation}
\begin{equation}\label{eq1.5} \sum_{m=1}^n \bm V_{mq}=n\bm V_q+o(n^{1-\epsilon})\; a.s., \;\; \text{ for all }\;
q=1,\ldots,d,
\end{equation}
where $\bm V_q=(V_{qkl})_{k,l=1}^d$, $q=1,\ldots,d$, are  $d\times
d$ non-negative definite matrices.
\end{assumption}

 By Assumption \ref{assumption1.1},   $\bm H$ has a   maximal
 eigenvalue $1$ and a corresponding   right eigenvector $\bm 1=(1,\ldots,1)$.
 Let $\lambda_2,\ldots,\lambda_d$ be other $d-1$ eigenvalues of $\bm H$.
 Then $\bm H$ has the following Jordan form
 decomposition
\begin{equation}\label{eq1.7}
 \bm T^{-1}\bm H \bm T =diag\left(
      1 , \bm J \right)\;\;\text{ and }\;\; \bm J=diag(\bm
      J_2,\ldots,\bm J_s)
 \end{equation}
 with
\begin{equation}\label{eq1.8}
\bm J_t= \left(\begin{matrix}
         \lambda_t &1           &  0     &\ldots     & 0 \\
         0         & \lambda_t  &1       &\ldots     & 0 \\
         \vdots    & \ldots     & \ddots &\ddots     & \vdots \\
         0         & 0          &\ldots  & \lambda_t & 1 \\
         0         & 0          &0       & \ldots    &\lambda_t
\end{matrix}\right),\end{equation}
where $\bm T=(\bm t_1^{\prime},\ldots,\bm t_d^{\prime})$ and $\bm
t_1=\bm 1$. Denote by
$\rho=\max\{Re(\lambda_2),\ldots,Re(\lambda_2)\}$, where
$Re(\lambda_k)$ is the real part of the complex number $\lambda_k$.
And denote the order of $\bm J_t$ by $\nu_t$ and $\nu=\max\{\nu_t:
Re(\lambda_t)=\rho\}$. Let $\bm v $ be the left eigenvector of $\bm
H$   associated with the positive maximal eigenvalue $1$ and satisfy
$\bm v \bm 1^{\prime}=1$. We denote $\widetilde{\bm H}= \bm H-\bm
1^{\prime}\bm v$, $\bm \Sigma_1=diag(\bm v)-\bm v^{\prime}\bm v$,
$\bm \Sigma_2=\sum_{q=1}^d v_q\bm V_q$ and $\bm \Sigma=\bm
H^{\prime}\bm \Sigma_1\bm H+ \bm \Sigma_2$. For a $d$-dimensional
Brownian motion $\{\bm W_t; t\ge 0\}$ with a co-variance $\bm
\Lambda$, we denote the solution of the equation:
$$  \bm S_t
=\bm W_t+ \int_0^t\frac{\bm S_s\widetilde{\bm H} }{s}\; ds, \quad
t\ge 0, \quad \bm S_0=\bm 0 \eqno Equ1 $$ by $\{\bm
S_t=\text{Solut}(\text{Equ1},\bm W_t); t\ge 0\}$. In the next
section, we will show that $\bm S_t$ is well defined if $\rho<1/2$,
and
$$ \bm S_t=\int_0^t \bm W_s\frac{\widetilde{\bm H}}{s}\Big(\frac ts\Big)^{\widetilde {\bm H}} ds+\bm W_t
=\int_0^t (d \bm W_s)\Big(\frac ts\Big)^{\widetilde {\bm H}}. $$
where, for any $t>0$ and any matrix $\bm M$, $t^{\bm M}$ is defined
to be
$$\exp\{\bm M\ln t\}:=\sum_{k=0}^{\infty}\frac{\bm M^k(\ln t)^k}{k!}.  $$
Also we denote the solution of the equation:
$$  \widehat{\bm S}_t
=\bm W_t-\bm W_1+ \int_1^t\frac{\widehat{\bm S}_s\widetilde{\bm H}
}{s}\; ds, \quad t>0, \quad \widehat{\bm S}_1=\bm 0 \eqno Equ2 $$ by
$\{\widehat{\bm S}_t=\text{Solut}(\text{Equ2},\bm W_t); t>0\}$. And
we will show that $\widehat{\bm S}_t$ is well defined  and
$$ \widehat{\bm S}_t=\int_1^t \bm W_s\frac{\widetilde{\bm H}}{s}\Big(\frac ts\Big)^{\widetilde {\bm H}} ds
+\bm W_t-\bm W_1t^{\widetilde{\bm H}}=\int_1^t (d \bm W_s)\Big(\frac
ts\Big)^{\widetilde {\bm H}}.
$$

Now let $\bm B_{t1}$ and $\bm B_{t2}$ be two independent
$d$-dimensional  standard Brownian motions. Define $\bm
G_{ti}=\text{Solut}(\text{Equ1}, \bm B_{ti}\bm\Sigma_i^{1/2})$ and
$\widehat{\bm G}_{ti}=\text{Solut}(\text{Equ2}, \bm
B_{ti}\bm\Sigma_i^{1/2})$, $i=1,2$. Let $\bm G_t=\bm G_{t1}\bm H+\bm
G_{t2}$ and $\widehat{\bm G}_t=\widehat{\bm G}_{t1}\bm
H+\widehat{\bm G}_{t2}$. Then $\bm G_t=\text{Solut}(\text{Equ1}, \bm
B_{t1}\bm\Sigma_1^{1/2}\bm H+\bm B_{t2}\bm\Sigma_2^{1/2})$ and
$\widehat{\bm G}_t=\text{Solut}(\text{Equ2}, \bm
B_{t1}\bm\Sigma_1^{1/2}\bm H+\bm B_{t2}\bm\Sigma_2^{1/2})$.

The next two theorems are on the strong approximation.
\begin{theorem}\label{theorem1.2} Suppose $\rho<1/2$.  Under Assumptions
\ref{assumption1.1} and \ref{assumption1.3}, there are two
independent $d-$dimensional standard Brownian motions $\bm B_{t1}$
and $\bm B_{t2}$ ({\em possibly in an  enlarged probability space
with  the process $(\bm Y_n, \bm N_n)$ being redefined without
changing the distribution}) such that for some $\gamma>0$,
\begin{eqnarray}\label{eq1.13}
 \bm Y_n-n\bm v &=&
 \bm G_{n1}\bm H+\bm G_{n2}+ o(n^{1/2-\tau\wedge\gamma}) \; a.s.,
 \\
\label{eq1.14} \bm N_n-n \bm v&=& \bm G_{n1}+\int_0^n \frac{\bm
G_{x2}}{x}\; dx(\bm I-\bm 1^{\prime}\bm
v)+o(n^{1/2-\tau\wedge\gamma}) \; a.s.,
\end{eqnarray}
where $\tau\wedge\gamma=\min\{\tau,\gamma\}$.
\end{theorem}

\begin{theorem}\label{theorem1.3} Suppose $\rho=1/2$ and Assumptions
\ref{assumption1.1} and
\ref{assumption1.3} are satisfied. Further, assume that
\begin{equation}\label{eq1.18}
\sum_{m=1}^{\infty}\frac{\| \bm H_m-\bm H\|}{m^{1/2}}<\infty.
\end{equation}
Then  there are two independent $d-$dimensional standard Brownian
motions $\bm B_{t1}$ and $\bm B_{t2}$ ({\em possibly in an  enlarged
probability space with  the process $(\bm Y_n, \bm N_n)$ being
redefined without changing the distribution}) such that
\begin{eqnarray}\label{eq1.15}
 \bm Y_n- n\bm v&=&\widehat{\bm G}_{n1} \bm H +\widehat{\bm G}_{n2}+O(n^{1/2}\log^{\nu-1}n) \; a.s.
 \\
\label{eq1.16} \bm N_n- n \bm v&=&\widehat{\bm G}_{n1}+\int_1^n
\frac{\widehat{\bm G}_{x2}}{x}\; dx(\bm I-\bm 1^{\prime}\bm
v)+O(n^{1/2}\log^{\nu-1}n) \; a.s.
\end{eqnarray}
Also,
\begin{equation}\label{eq1.17} \Big(\widehat{\bm G}_{t1}+
\int_1^t \frac{\widehat{\bm G}_{x2}}{x}\; dx(\bm I-\bm 1^{\prime}\bm
v)\Big)\bm H =\widehat{\bm G}_{t1}\bm H+\widehat{\bm G}_{t2} -\bm
B_{t2}\bm\Sigma_2^{1/2}. \end{equation}
 \end{theorem}

\begin{remark} The condition (\ref{eq1.18}) is used by Bai and Hu (2005) to obtain the asymptotic normality.
It is easily seen that it implies the condition
(\ref{eqassumption1.1.1}) with $\tau=0$.  Bai and Hu (2005) also
assumed that $\bm H_m\bm 1^{\prime}= \bm 1^{\prime}$ and $\bm H_m\to
\bm H$.

\end{remark}
\begin{remark} By (\ref{eq1.15})-(\ref{eq1.17}), under the
assumptions in Theorem \ref{theorem1.3},
$$ \bm Y_n-n \bm v=(\bm N_n-n \bm v)\bm H +O(n^{1/2}\log^{\nu-1} n)\;\; a.s.\;\; \text{ if } \nu>1. $$
\end{remark}
The proof of Theorems \ref{theorem1.2} and \ref{theorem1.3} will be
given in the last section. Before that, we give some properties of
the limit processes and several application of these approximations.


\section{Properties of the limit processes.} \setcounter{equation}{0}

This section will give several properties of  the solutions of
equations (Equ1) and (Equ2) and Gaussian processes $\bm G_t$, $\bm
S_t$, etc. By combining  these properties with the approximation
theorems in the above section we can obtain important properties of
the urn models which will be given in the next section. The
properties of the Gaussian processes will be also used in the proofs
of the approximation theorems.

We
first need two lemmas.

\begin{lemma}\label{lemma2.7} Under Assumption \ref{assumption1.1}, there
exists a constant $C$ such that for any $a\ge 1$,
\begin{equation}\label{eq2.10}\| a^{\widetilde{\bm H}}\|\le C
a^{\rho}\log^{\nu-1}a. \end{equation}
\end{lemma}
\noindent{\bf Proof} It is obvious that $\widetilde{\bm H} =\bm
Tdiag(0, \bm J)\bm T^{-1}: =\bm T \widetilde{\bm J}\bm T^{-1}$. It
follows that
$$a^{\widetilde{\bm H}}
=\bm T
 a^{\widetilde{\bm J}}
\bm T^{-1}.$$ So it is enough to show that for any $a>1$, $ \|
a^{\bm J_t}\|\le C a^{Re(\lambda_t)}\log^{\nu_t-1}a$.   Denote $\bm
J_t=\lambda_t \bm I +\overline{\bm I}_t $ where
\begin{equation}\label{eq2.11}\overline{\bm I}_t= \left(\begin{matrix}
         0         &1           &  0     &\ldots     & 0 \\
         0         &  0          &1       &\ldots     & 0 \\
         \vdots    & \ldots     & \ddots &\ddots     & \vdots \\
         0         & 0          &\ldots  &   0       & 1 \\
         0         & 0          &0       & \ldots    & 0
\end{matrix}\right).\end{equation}
Then $\| a^{\bm J_t}\|=\|a^{\lambda_s} a^{\overline{\bm I}_t}\| \le
C a^{Re(\lambda_t)}\|a^{\overline{\bm I}_t}\|$.  Obviously,
$$\overline{\bm I}_t^2=\left(\begin{matrix}
         0         & 0           &  1     &\ldots     & 0 \\
         0         &  0          &  0     &\ldots     & 0 \\
         \vdots    & \vdots      & \vdots &\dots     & \vdots \\
         0         & 0          & 0       & \ldots    & 1 \\
         0         & 0          & 0       &  \ldots    & 0 \\
         0         & 0          &0       & \ldots    & 0
\end{matrix}\right),  \ldots,
\overline{\bm I}_t^{\nu_s-1}=\left(\begin{matrix}
         0         & 0           &  0     &\ldots     & 1 \\
         0         &  0          &  0     &\ldots     & 0 \\
         \vdots    & \vdots      & \vdots &\dots     & \vdots \\
         0         & 0          & 0       & \ldots    & 0 \\
         0         & 0          & 0       &  \ldots    & 0 \\
         0         & 0          &0       & \ldots    & 0
\end{matrix}\right)
$$
and $\overline{\bm I}_t^{\nu_s}=\bm 0$. It follows that
$$ a^{\overline{\bm I}_t}=\sum_{k=0}^{\infty}
\frac{\overline{\bm I}_t^k(\ln a)^k}{k!} =\sum_{k=0}^{\nu_t-1}
\frac{\overline{\bm I}_t^k(\ln a)^k}{k!}. $$ Then
$$ \|a^{\overline{\bm I}_t}\|\le
C\sum_{k=0}^{\nu_t-1} \frac{(\ln a)^k}{k!}\le C(\log a)^{\nu_s-1}.
$$ Hence (\ref{eq2.10}) is proved.

\begin{lemma}\label{lemma3.1} For any $a\ge 0$, the equation
$$\bm Z_t=\int_a^t\frac{\bm Z_s}{s}\widetilde{\bm H} ds, \quad \bm Z_a=\bm 0 $$
or equivalently \begin{equation} \label{eqlem3.1.1}d\bm
Z_t=\frac{\bm Z_t}{t}\widetilde{\bm H} dt, \quad \bm Z_a=\bm 0
\end{equation} has an
unique solution $\bm Z_t \equiv \bm 0$.
\end{lemma}
\noindent{\bf Proof} It is obvious that  $\widetilde{\bm H}=\bm
H-\bm 1^{\prime}\bm v$ has the Jordan form decomposition
$$ \bm T^{-1} \widetilde{\bm H} \bm T
=diag(0,\bm J_2,\ldots,\bm J_s) $$
and (\ref{eqlem3.1.1}) is equivalent to
$$ \bm
Z_t \bm T =\frac{\bm Z_t \bm T}{t} diag(0,\bm J_2,\ldots,\bm J_s) dt, \quad \bm Z_a\bm T=\bm 0. $$
On the other hand,  for each $s$,
$$ d \widetilde{\bm Z}_t^{(s)} = \frac{\widetilde{\bm Z}_t^{(s)}}{t} \bm J_{s} dt,
\quad \widetilde{\bm Z}_a^{(s)}=\bm 0 $$ has an unique solution $
\widetilde{\bm Z}_t^{(s)} \equiv \bm 0$. The proof is completed.

  From this Lemma, it follows
that the solutions of (Equ1) and (Equ2) are unique. The following
two propositions tells us that the solutions exist.

\begin{proposition}\label{lemma3.2} Let $\{\bm W_t; t\ge 0\}$ be a
$d$-dimensional Brownian motion with a co-variance matrix $\bm
\Lambda$. Suppose that Assumption \ref{assumption1.1} is satisfied
and $\rho<1/2$. Then the unique solution   $\bm S_t
=\text{Solut}(\text{Equ1}, \bm W_t)$  of the equation (Equ1) is
\begin{equation}\label{eq3.1}
\bm S_t=\int_0^t \bm W_x \frac{\widetilde{\bm H}}{x} \big(\frac tx
\big)^{\widetilde{\bm H}} dx +\bm W_t.  \end{equation}
 Also
\begin{equation}\label{eq3.2}\bm S_t=\Big(\int_0^t (d\bm W_x)
x^{-\widetilde{\bm H} } \Big) t^{\widetilde{\bm H} }\quad a.s.
\end{equation} Furthermore, with probability one $\bm S_t$ is
continuous on $[0,\infty)$.
\end{proposition}

\noindent{\bf Proof} Fist, since $\|x^{-\widetilde{\bm H}}\|\le C
x^{-\rho}(\log x^{-1})^{\nu -1}$ for $0<x<1$ by Lemma
\ref{lemma2.7}, and $\bm W_x \overset{a.s.}= O(\sqrt{x\log\log
x^{-1} })$ as $x\to 0$,  we have
$$\bm W_x x^{-\widetilde{\bm H}}=O(1)x^{1/2-\rho}(\log x^{-1})^{\nu }
\to 0 \quad a.s., $$
$$
\bm W_x  \frac {\widetilde{\bm H}}{x} \big(\frac
tx\big)^{\widetilde{\bm H}} = O(1)x^{-1/2-\rho}(\log x^{-1})^{\nu}
\quad a.s. $$
 as  $x\to 0$. Since $-1/2-\rho>-1$ and
$\bm W_x \frac {\widetilde{\bm H}}{x} \big(\frac
tx\big)^{\widetilde{\bm H}}$ is continuous on $(0,\infty)$, it
follows that the integral $\int_0^t \bm W_x \frac{\widetilde{\bm
H}}{x}\big(\frac tx\big)^{\widetilde{\bm H}}dx $ exists, and then
$\bm S_t$ in (\ref{eq3.1}) is well defined and
\begin{eqnarray*}
&&\int_0^t \bm W_x \frac{ \widetilde{\bm H} }{x}\big(\frac tx\big)^{
\widetilde{\bm H} } dx =\Big(\int_0^t \bm W_x d(x^{ -\widetilde{\bm
H} } )\Big)
  t^{ \widetilde{\bm H} } \\
&=& -\bm W_x x^{ -\widetilde{\bm H} }\big|_0^t t^{ \widetilde{\bm
H}} +\Big(\int_0^t d(\bm W_x) x^{-\widetilde{\bm
H}}\Big)t^{\widetilde{\bm H}} =-\bm W_t+\Big( \int_0^t d(\bm W_x)
x^{-\widetilde{\bm H} }\Big) t^{ \widetilde{\bm H} }.
\end{eqnarray*}
It follows that (\ref{eq3.2}) is true. Now we show that $\bm S_t$ is
the solution of equation (Equ1). Note that
$$\bm S_t=O(1)\int_0^t(x\log\log x^{-1})^{1/2}(t/x)^{\rho}(\log(t/x))^{\nu -1}dx
=O(1)t^{1/2-\rho}(\log t^{-1})^{\nu} $$
as $t\to \infty$. It follows
that $\bm S_0=0$, the integral $\int_0^t\frac{\bm S_s}{s}ds$ exists
and
\begin{eqnarray*}
&&\int_0^t\frac{\bm S_s}{s}ds = \int_0^t\frac{ds}{s} \int_0^s \bm
W_x\frac{\widetilde{\bm H}}{x}\big(\frac sx\big)^{\widetilde{\bm
H}}dx
+\int_0^t \frac{\bm W_s}{s}ds \\
&=& \int_0^t \bm W_x x^{-\widetilde{\bm H}-1} dx
 \int_x^t \widetilde{\bm H} s^{\widetilde{\bm H}-1}ds
+\int_0^t \frac{\bm W_s}{s}ds \\
&=& \int_0^t \bm W_x x^{-\widetilde{\bm H}-1} dx
  s^{\widetilde{\bm H}}\big|_x^t
+\int_0^t \frac{\bm W_s}{s}ds \\
&=& \int_0^t \bm W_x x^{-\widetilde{\bm H}-1} dx
  t^{\widetilde{\bm H}}-\int_0^t \frac{\bm W_x}{x}dx
+\int_0^t \frac{\bm W_s}{s}dx =\int_0^t \bm W_x\frac 1x \big(\frac
tx\big)^{\widetilde{\bm H}}dx.
\end{eqnarray*}
Then
$$\int_0^t\frac{\bm S_s\widetilde{\bm H} }{s}ds
=\int_0^t \bm W_x\frac{\widetilde{\bm H}}{x}
  \big(\frac tx\big)^{\widetilde{\bm H}}dx
=\bm S_t-\bm W_t. $$ So, $\bm S_t$ is the solution of equation
(Equ1). Finally, the continuity of $\bm S_t$ follows from the
continuity of the Brownian motion $\bm W_t$.

\begin{proposition}\label{lemma3.3} Let $\{\bm W_t; t\ge 0\}$ be a
$d$-dimensional Brownian motion with some co-variance matrix.
Suppose Assumption \ref{assumption1.1} is satisfied. Then the unique
solution $\widehat{\bm S}_t =\text{Solut}(\text{Equ2}, \bm W_t)$ of
the equation (Equ2)  is
\begin{eqnarray}\label{eq3.3}\widehat{\bm S}_t=\int_1^t \bm W_x \frac{\widetilde{\bm H}}{x} \big(\frac tx \big)^{\widetilde{\bm H}} dx +\bm W_t-\bm W_1t^{\widetilde{\bm H}}. \end{eqnarray}
Also
\begin{eqnarray}\label{eq3.4}\widehat{\bm S}_t=\Big(\int_1^t (d\bm W_x) x^{-\widetilde{\bm H} } \Big)
t^{\widetilde{\bm H} }\quad a.s. \end{eqnarray}
 Furthermore, with
probability one $\widehat{\bm S}_t$ is continuous on $(0,\infty)$.
\end{proposition}

\noindent{\bf Proof} The proof is similar to Proposition \ref{lemma3.2} and so   omitted.

\begin{proposition} \label{proposition2.3}Let $\{\bm W_t; t\ge 0\}$ be a
$d$-dimensional Brownian motion with some co-variance matrix.
Suppose Assumption \ref{assumption1.1} is satisfied. If $\rho<1/2$,
then for $\bm S_t =\text{Solut}(\text{Equ1}, \bm W_t)$ we have
\begin{equation}\label{eqproposition2.3.1}
 \int_0^n \frac{\bm S_t}t dt =\sum_{m=1}^{n-1}\frac{\bm S_m}{m}+O(1) \;\; a.s.
 \end{equation}
If $\rho<1$, then for $\widehat{\bm S}_t =\text{Solut}(\text{Equ2},
\bm W_t)$ we have
\begin{equation}\label{eqproposition2.3.2}
 \int_1^n \frac{\widehat{\bm S}_t}t dt =\sum_{m=1}^{n-1}\frac{\widehat{\bm S}_m}{m}+O(1) \;\; a.s.
 \end{equation}
\end{proposition}

\noindent{\bf Proof} We only give  a proof of
(\ref{eqproposition2.3.2}) since the proof of
(\ref{eqproposition2.3.1}) is similar. First, from (\ref{eq3.4}) it
follows that for all $t>1$,
\begin{align*}
 \|\Var( \widehat{\bm S}_t)\|=& \left\| \int_1^t \big(\frac{t}{x}\big)^{\widetilde{\bm H}^{\prime}}\Var(\bm W_1)
 \big(\frac{t}{x}\big)^{\widetilde{\bm H}}dx\right\|
 \le    C \int_1^t (t/x)^{2\rho}\log^{2\nu-2}(t/x) dx \\
 \le & \left\{\begin{matrix}C t \log^{2\nu -1} t , & \text{ if } \rho=1/2\\
 Ct, & \text{ if } \rho<1/2  \\
 C t^{2\rho} \log^{2\nu -2}\log t, & \text{ if } \rho>1/2
 \end{matrix}\right\}
 \quad \le C t^{2(\rho\vee
 \frac{1}{2})}\log^{2\nu-1} t.
 \end{align*}
 So, $\ep\|\widehat{\bm S}_t\|\le C t^{\rho\vee\frac{1}{2}} \log^{\nu -1/2} t$. According to equation (Equ2),
 $$ \widehat{\bm S}_t-\widehat{\bm S}_s = \bm W_t-\bm W_s+\int _s^t \frac{\widehat{\bm S}_x\widetilde{\bm H}}{x} dx, \;\; t\ge s\ge 1. $$
It follows that
\begin{align*}
&\sum_{m=1}^{\infty} \int_m^{m+1}\Big (\frac{\widehat{\bm S}_t}t -\frac{\widehat{\bm S}_m}m\Big) dt \\
=& \sum_{m=1}^{\infty}\int_m^{m+1}\widehat{\bm S}_t\big(\frac{1}{t}-\frac{1}{m}\big)dt+\sum_{m=1}^{\infty} \int_m^{m+1}\frac{\widehat{\bm S}_t-\widehat{\bm S}_m}{m}dt\\
=& \sum_{m=1}^{\infty}\int_m^{m+1}\widehat{\bm S}_t\big(\frac{1}{t}-\frac{1}{m}\big)dt
+\sum_{m=1}^{\infty} \int_m^{m+1}\frac{\bm W_t-\bm W_m}{m}dt\\
&+\sum_{m=1}^{\infty} \frac{1}{m} \int_m^{m+1}\int_m^t\frac{\widehat{\bm S_x}\widetilde{\bm H}}{x}dxdt.
\end{align*}
The first  and the third term above are a.s. convergent because
$$\sum_{m=1}^{\infty}\int_m^{m+1}\ep\|\widehat{\bm S}_t\|\big|\frac{1}{t}-\frac{1}{m}\big|dt\le C \sum_{m=1}^{\infty}
\frac{(m+1)^{\rho\vee \frac{1}{2}}\log^{v-1/2}(m+1)}{m^2}<\infty, $$
and
$$\sum_{m=1}^{\infty} \frac{1}{m} \int_m^{m+1}\int_m^t\frac{\ep\|\widehat{\bm S_x}\widetilde{\bm H}\|}{x}dxdt
\le C \sum_{m=1}^{\infty} \frac{(m+1)^{\rho\vee
\frac{1}{2}}\log^{v-1/2}(m+1)}{m^2}<\infty.
$$ The second term is a.s. convergent because it is an infinite
series of independent normal random variables with
$$\sum_{m=1}^{\infty} \left\|\Var\left\{\int_m^{m+1}\frac{\bm W_t-\bm W_m}{m}dt\right\}\right\|\le C\sum_{m=1}^{\infty} \frac{1}{m^2}<\infty. $$
It follows that
$$\sum_{m=1}^{n-1} \int_m^{m+1}\Big (\frac{\widehat{\bm S}_t}t -\frac{\widehat{\bm S}_m}m\Big) dt=O(1) \;\; a.s. $$
The proof of (\ref{eqproposition2.3.2}) is completed.

\bigskip

Propositions \ref{lemma3.2} and \ref{lemma3.3} give the solutions of
equations (Equ1) and (Equ2). To give further   properties of the
$\bm G_t$ and $\widehat{\bm G}_t$, we need the analytic
representation of the solutions. Recall $\bm T=(\bm
t_1^{\prime},\ldots, \bm t_d^{\prime})$, where $\bm t_1=\bm 1$. Let
$\{\bm W_t; t\ge 0\}$ be a $d$-dimensional Brownian motion with some
co-variance matrix $\bm \Lambda$.  First we consider (Equ1). Let
$\{\bm S_t=\text{Solut}(\text{Equ1},\bm W_t); t\ge 0\}$ be the
solution of (Equ1) and $\bm U_t=\bm S_t \bm T$. Then $\bm U_t$ is
the unique solution of the equation
\begin{eqnarray}\label{eq3.17}  \bm U_t
=\bm W_t\bm T+ \int_0^t\frac{\bm U_s\widetilde{\bm J} }{s}\; ds
\quad t\ge 0, \quad \bm U_0=\bm 0 \end{eqnarray}
 Note that
$\widetilde{\bm J}=  diag(0, \bm J)= diag(0,\bm J_2,\ldots,\bm J_s),
$ where $\bm J_i$'s are defined as in (\ref{eq1.8}). Write $\bm
U_t=(U_{t1}, \bm U_t^{(2)},\ldots,\bm U_t^{(s)})$, where $\bm
U_t^{(i)}=(U_{t1}^{(i)},\ldots, U_{t\nu_i}^{(i)})$ is the vector
which contains $\nu_i$ coordinate variables corresponding to $\bm
J_i$. Also write $\bm T=(\bm 1^{\prime},\bm T^{(2)},$  $\ldots,$
$\bm T^{(s)})$, where $\bm T^{(i)}=(\bm t_{i1}^{\prime},\ldots,\bm
t_{i\nu_i}^{\prime})$ is the $\nu_i\times d$ matrix which contains
$\nu_i$ columns of $\bm T$ corresponding to $\bm J_i$. Obviously,
$U_{tj}^{(i)}=U_{t,1+\nu_2+\ldots+\nu_{i-1}+j}$ and $\bm t_{ij}=\bm
t_{1+\nu_2+\ldots+\nu_{i-1}+j}$. It is easily seen that
(\ref{eq3.17}) is equivalent to
\begin{equation}\label{eq3.18}\begin{array}{rl}
&U_{t1}=\bm W_t\bm 1^{\prime} \\
& d \; U_{t1}^{(i)}=d(\bm W_t \bm t_{i1}^{\prime})+\lambda_i \frac
{U_{t1}^{(i)}}{t}\; dt,
\quad  U_{01}^{(i)}=0,\\
& d \; U_{tj}^{(i)}=d(\bm W_t \bm t_{ij}^{\prime})+\frac
{U_{t,j-1}^{(i)}}{t}+\lambda_i \frac {U_{tj}^{(i)}}{t}\; dt,
\quad  U_{0j}^{(i)}=0,\\
& j=2,\ldots, \nu_i; \quad i=2,\ldots,s.
\end{array} \end{equation}
On can show that the solution of equation (\ref{eq3.18}) is
\begin{equation}\label{eq3.19}\begin{array}{rl}
&U_{t1}=\bm W_t\bm 1^{\prime} \\
U_{t1}^{(i)}= & t^{\lambda_i} \int_0^t \frac{d(\bm W_x \bm t_{i1}^{\prime})}{x^{\lambda_i}},\\
U_{tj}^{(i)}= & t^{\lambda_i}\int_0^t \frac{d(\bm W_x \bm
t_{ij}^{\prime})}{x^{\lambda_i}},
             +t^{\lambda_i} \int_0^t \frac{U_{x,j-1}^{(i)}}{x^{1+\lambda_i}}\;dx,    \\
              & j=2,\ldots, \nu_i; \quad i=2,\cdots,s.
\end{array} \end{equation}
Putting all the $U$'s to $\bm S_t=\bm U_t\bm T^{-1}$, we obtain the
solution of (Equ1).

Similarly, we have $\widehat{\bm S}_t=\widehat{\bm U}_t \bm T^{-1}$, where
\begin{equation}\label{eq3.20}  \bm U_t =\bm W_t\bm T-\bm W_1 \bm T +
\int_1^t\frac{\widehat{\bm U}_s\widetilde{\bm J} }{s}\; ds \quad t>
0, \quad \widehat{\bm U}_1=\bm 0 \end{equation}

and, $\widehat{\bm U}_t=(\widehat U_{t1}, \widehat{\bm
U}_t^{(2)},\ldots,\widehat{\bm U}_t^{(s)})$, $\widehat{\bm
U}_t^{(i)}=(\widehat U_{t1}^{(i)},\ldots, \widehat
U_{t\nu_i}^{(i)})$,
\begin{equation}\label{eq3.21}\begin{array}{rl}
\widehat U_{t1}= &(\bm W_t-\bm W_1)\bm 1^{\prime}, \\
\widehat U_{t1}^{(i)}= & t^{\lambda_i} \int_1^t \frac{d(\bm W_x \bm t_{i1}^{\prime})}{x^{\lambda_i}},\\
\widehat U_{tj}^{(i)}= & t^{\lambda_i}\int_1^t \frac{d(\bm W_x \bm
t_{ij}^{\prime})}{x^{\lambda_i}},
             +t^{\lambda_i} \int_1^t \frac{U_{x,j-1}^{(i)}}{x^{1+\lambda_i}}\;dx,    \\
              & j=2,\ldots, \nu_i; \quad i=2,\cdots,s.
\end{array} \end{equation}

\begin{proposition}\label{proposition2.4}
Under Assumption \ref{assumption1.1} and $\rho<1/2$,
\begin{equation}\label{eqproposition2.4.1}  \Var\big\{\big( \bm G_{t1}\bm H+ \bm
G_{t2},  \bm G_{t1}+\int_0^t \frac{ \bm G_{x2}
}{x}dx (\bm I-\bm 1^{\prime}\bm v)\big)\big\}=t\bm\Gamma
\end{equation}
with
\begin{align}\label{eqproposition2.4.2}
\bm \Gamma=
\Var\big\{\big( & \bm G_{11}\bm H+ \bm
G_{12},  \bm G_{11}+\int_0^1 \frac{ \bm G_{x2}
}{x}dx (\bm I-\bm 1^{\prime}\bm v)\big)\big\} \\
&=:\left(\begin{matrix} \bm \Gamma^{(11)} &  \bm \Gamma^{(12)} \\
          \bm \Gamma^{(21)} & \bm \Gamma^{(22)}\end{matrix}\right)\nonumber
\end{align}
and
\begin{align*}
\bm \Gamma^{(11)}=& \int_0^1 \left(\frac{1}{x}\right)^{\widetilde{\bm H}^{\prime}} \big(\bm H^{\prime}\bm\Sigma_1\bm H +\bm\Sigma_2\big)\left(\frac{1}{x}\right)^{\widetilde{\bm H} }dx,\\
\bm \Gamma^{(22)}=& \int_0^1 \left(\frac{1}{x}\right)^{\widetilde{\bm H}^{\prime}}  \bm\Sigma_1 \left(\frac{1}{x}\right)^{\widetilde{\bm H} }dx\\
& +(\bm I-\bm 1^{\prime}\bm v)^{\prime}\int_0^1\left[\int_x^1
\frac{1}{y} \left(\frac{y}{x}\right)^{\widetilde{\bm
H}}dy\right]^{\prime}\bm\Sigma_2
\left[\int_x^1 \frac{1}{y} \left(\frac{y}{x}\right)^{\widetilde{\bm H}}dy\right]dx(\bm I-\bm 1^{\prime}\bm v),\\
\bm \Gamma^{(12)}=&\overline{\bm \Gamma^{(12)}} =\bm H^{\prime}
\int_0^1 \left(\frac{1}{x}\right)^{\widetilde{\bm H}^{\prime}}
\bm\Sigma_1 \left(\frac{1}{x}\right)^{\widetilde{\bm H} }dx\\
 &\qquad \quad +(\bm I-\widetilde{\bm H}^{\prime})^{-1}\int_0^1 \left(\frac{1}{x}\right)^{\widetilde{\bm H}^{\prime}}
  \bm\Sigma_2 \left(\frac{1}{x}\right)^{\widetilde{\bm H} }dx (\bm I-\bm 1^{\prime}\bm v).
\end{align*}
\end{proposition}

\noindent {\bf Proof} Let $\{\bm W_t; t\ge 0\}$ be a $d$-dimensional
Brownian motion with some co-variance matrix $\bm \Lambda$, $\{\bm
S_t=\text{Solut}(\text{Equ1},\bm W_t); t\ge 0\}$ be the solution of
(Equ1). Notice $\{T^{-1/2}\bm W_{Tt}, t\ge 0\}$ and $\{\bm W_t,t\ge
0\}$ are identical distributed.
 So  $\{T^{-1/2}\bm S_{Tt}, t\ge 0\}$ and $\{\bm S_t,t\ge 0\}$ are identical distributed. Hence
 (\ref{eqproposition2.4.1}) is true.
  By (\ref{eq3.2}),
 \begin{align*}
\int_0^t\frac{\bm S_y}{x}dy=& \int_0^t\left[ \frac{1}{y}\int_0^y d\bm W_x\left(\frac{y}{x}\right)^{\widetilde{\bm H}}\right]dy
= \int_0^t d\bm W_x\left[\int_x^t \frac{1}{y}\left(\frac{y}{x}\right)^{\widetilde{\bm H}}dy\right].
\end{align*}
It follows that
$$
  \Var\{\bm S_1\}= \int_0^1 x^{-\widetilde{\bm H}^{\prime} }\bm\Lambda
x^{-\widetilde{\bm H} }dx,  $$
$$\Cov\left\{\bm S_t,\bm S_s\right\}=
\int_0^s \left(\frac{t}{x}\right)^{\widetilde{\bm H}^{\prime}}\bm\Lambda
\left(\frac{s}{x}\right)^{\widetilde{\bm H}}dx
=
s\left(\frac{t}{s}\right)^{\widetilde{\bm H}^{\prime}}\Var\{\bm S_1\}, \;\; t\ge s, $$
$$\Var\left\{\int_0^1\frac{\bm S_y}{x}dy\right\}
=\int_0^1\left[\int_x^1 \frac{1}{y} \left(\frac{y}{x}\right)^{\widetilde{\bm H}}dy\right]^{\prime}\bm\Lambda
\left[\int_x^1 \frac{1}{y} \left(\frac{y}{x}\right)^{\widetilde{\bm H}}dy\right]dx,
$$
\begin{align*}
\Cov\left\{ \bm S_1, \int_0^1\frac{\bm S_y}{x}dy\right\}=& \int_0^1\frac{\Cov\{\bm S_1,\bm S_y\}}{y} dy
=\int_0^1\left(\frac{1}{y}\right)^{\widetilde{\bm H}^{\prime}}dy \Var\{\bm S_1\}\\
=&(\bm I-\widetilde{\bm H}^{\prime})^{-1}\Var\{\bm S_1\}.
\end{align*}
The proof is now completed by noticing  the independence of $\bm
G_{t1}$ and $\bm G_{t2}$.

\begin{proposition} \label{proposition2.5}
Under Assumption \ref{assumption1.1} and $\rho=1/2$, the limit
\begin{equation}\label{eq7.2}\widetilde{\bm \Gamma}=\lim_{t\to \infty} t^{-1}(\log
t)^{1-2\nu} \Var\big\{\big(\widehat{\bm G}_{t1}\bm H+\widehat{\bm
G}_{t2}, \widehat{\bm G}_{t1}+\int_1^t \frac{\widehat{\bm G}_{x2}
}{x}dx (\bm I-\bm 1^{\prime}\bm v)\big)\big\}
\end{equation}
exists, and
\begin{equation}\label{eq7.4}\widetilde{\bm \Gamma}
=:\left(\begin{matrix} \widetilde{\bm \Gamma}^{(11)} &  \widetilde{\bm \Gamma}^{(12)} \\
           \widetilde{\bm \Gamma}^{(21)} & \widetilde{\bm \Gamma}^{(22)}\end{matrix}\right),
\end{equation}
where
\begin{equation}\label{eq7.3}
\begin{array}{rl}
 (\bm T^{\ast}\widetilde{\bm \Gamma}^{(11)}\bm T)_{ij} & =
\frac 1{((\nu-1)!)^2} \frac 1{2\nu-1}
    \big(|\lambda_l|^2\overline{\bm t_{i1}}\bm\Sigma_1\bm t_{j1}^{\prime}
     +\overline{\bm t_{i1}}\bm\Sigma_2\bm t_{j1}^{\prime} \big), \\
 (\bm T^{\ast}\widetilde{\bm \Gamma}^{(22)}\bm T)_{ij}& =
\frac 1{((\nu-1)!)^2} \frac 1{2\nu-1}
    \big( \overline{\bm t_{i1}}\bm \Sigma_1\bm t_{j1}^{\prime}
     +|\lambda_l|^{-2}\overline{\bm t_{i1}}\bm \Sigma_2\bm t_{j1}^{\prime}\big), \\
 (\bm T^{\ast}\widetilde{\bm \Gamma}^{(12)}\bm T)_{ij} & =
\overline{  (\bm T^{\ast}\widetilde{\bm \Gamma}^{(21)}\bm T)_{ij} }
=\frac 1{((\nu-1)!)^2} \frac 1{2\nu-1}
    \big( \overline{\lambda}_l\overline{\bm t_{i1}}\bm\Sigma_1\bm t_{j1}^{\prime}
     +\lambda_l^{-1}\overline{\bm t_{i1}}\bm\Sigma_2\bm t_{j1}^{\prime}\big)
\end{array}\end{equation}
whenever $i=j=1+\nu_2+\ldots+\nu_l$ and $Re(\lambda_l)=1/2$,
$\nu_l=\nu$, and $ (\bm T^{\ast}\widetilde{\bm \Gamma}^{(uv)}\bm
T)_{ij}=0$ for otherwise $u,v=1,2$. Here $\overline{\bm a}$ is the
conjugate vector of a complex   vector $\bm a$.

\end{proposition}

\noindent{\bf Proof}  Let $\bm W_t$ be a $d$-dimensional Brownian
motion with a co-variance $\bm \Lambda$, $\widehat{\bm S}_t$ a
solution of (Equ2) and $\widehat{\bm U}_t=\widehat{\bm S}_t\bm T$.
Then by (\ref{eq3.20}) and Proposition \ref{lemma3.3},
$$\widehat U_{t1}=(\bm W_t-\bm W_1)\bm 1^{\prime}\overset{L_2}
=o( t^{1/2}\log^{v-1/2}t )$$
and
\begin{eqnarray*}
\widehat{\bm U}_t^{(i)}&=& \bm W_t\bm T^{(i)}-\bm W_1\bm T^{(i)}
t^{\bm J_i}
+\int_1^t \bm W_x\bm T^{(i)}\frac 1x \big(\frac tx\big)^{\bm J_i} \bm J_i \; dx \\
 &\overset{L_2}=& o(t^{1/2}\log^{v-1/2} t)
+\sum_{k=0}^{\nu_i-1}\frac 1{k!}
 \int_1^t \bm W_x\bm T^{(i)}\frac 1x \big(\frac tx\big)^{\lambda_i} \log^k \frac tx \overline{\bm I}_i^k\bm J_i \;
 dx
 \\
&=& o(t^{1/2}\log^{v-1/2} t) +\sum_{k=0}^{\nu_i-1}\frac
1{k!}\lambda_i
 \int_1^t \bm W_x\bm T^{(i)}\frac 1x \big(\frac tx\big)^{\lambda_i} \log^k \frac tx \overline{\bm I}_i^k\;
 dx,
\end{eqnarray*}
where $\overline{\bm I}_i$ is defined as in (\ref{eq2.11}). It is
easily seen that
\begin{eqnarray*} \int_1^t \bm W_x\bm T^{(i)}\frac 1x \big(\frac tx\big)^{\lambda_i}\log^k \frac tx
&\overset{L_2}=& O(1)\int_1^tx^{1/2}\frac 1x \big(\frac tx\big)^{Re(\lambda_i)}\log^k \frac tx \\
&=& \begin{cases} O(t^{1/2}), & \text{ if } Re(\lambda_i)<1/2, \\
         O(t^{1/2}\log^{k+1}t ), & \text{ if } Re(\lambda_i)=1/2. \end{cases}
\end{eqnarray*}
So
$$\widehat{\bm U}_t^{(i)}\overset{L_2}=
\frac{\lambda_i}{(\nu_i-1)!}
 \int_1^t \bm W_x\bm T^{(i)}\frac 1x \big(\frac tx\big)^{\lambda_i} \big(\log^{\nu_i-1}
 \frac tx \big)\overline{\bm I}_i^{\nu_i-1}\; dx
 +o(t^{1/2}\log^{v-1/2} t).
 $$
It follows that
\begin{eqnarray*}
&&\widehat U_{tj}^{(i)} \overset{L_2}= \begin{cases}
\frac{\lambda_i}{(\nu-1)!}
       \int_1^t \bm W_x\bm t_{i1}^{\prime}\frac 1x \big(\frac tx\big)^{\lambda_i}
       \log^{\nu-1} \frac tx \; dx, & \text{ if } Re(\lambda_i)=1/2,
         j= \nu_i=\nu,  \\
       0,&  \text{ otherwise}.
\end{cases} \\
&& \qquad \qquad +o(t^{1/2}\log^{v-1/2} t).
\end{eqnarray*}
Similarly,
\begin{eqnarray*}
\int_1^t\frac{\widehat U_{xj}^{(i)}}x dx &\overset{L_2}=&
\begin{cases} \frac{1}{(\nu-1)!}
       \int_1^t \bm W_x\bm t_{i1}^{\prime}\frac 1x \big(\frac tx\big)^{\lambda_i}
       \log^{\nu-1} \frac tx \; dx, & \text{ if } Re(\lambda_i)=1/2,
        j= \nu_i=\nu,  \\
       0,&  \text{ otherwise}.
\end{cases}\\
&&+o(t^{1/2}\log^{v-1/2} t)
\end{eqnarray*}
On the other hand, if
$Re(\lambda_i)=Re(\lambda_j)=1/2$ and $\nu_i= \nu_j=\nu$, then
\begin{eqnarray*}
&&\Cov\Big\{  \int_1^t \bm W_x\bm t_{i1}^{\prime}\frac 1x \big(\frac
tx\big)^{\lambda_i}
       \log^{\nu-1} \frac tx \; dx,
       \int_1^t \bm W_x\bm t_{j1}^{\prime}\frac 1x \big(\frac tx\big)^{\lambda_j}
       \log^{\nu-1} \frac tx \; dx\Big\} \\
&=& \begin{cases}
    (2\nu-1)^{-1}|\lambda_i|^{-2}\overline{\bm t_{i1}}\bm \Lambda \bm t_{j1}^{\prime}
 \big(1+o(1)\big) t\log^{2\nu-1} t& \text{ if } \lambda_i=\lambda_j, \\
\big( \frac 1{\overline{\lambda}_i}+\frac 1{\lambda_j})\frac
1{1-\overline{\lambda}_i-\lambda_j} \overline{\bm t_{i1}} \bm
\Lambda \bm t_{j1}^{\prime} \big(1+o(1)\big)
t^{2-\overline{\lambda}_i-\lambda_j}\log^{2\nu-2}t, & \text{ if }
\lambda_i\ne\lambda_j.\end{cases}
\end{eqnarray*}
It follows that
\begin{align*}
\lim_{t\to \infty} t^{-1}(\log t)^{1-2\nu}
& \Var\Big\{(\widehat{\bm G}_{t1}\bm H+\widehat{\bm G}_{t2})\bm T,
         \big(\widehat{\bm G}_{t1}
   +\int_1^t\frac {\widehat{\bm G}_{x1}}x  dx (\bm I-\bm 1^{\prime}\bm v)\big)\bm T
\Big\}\\
=:&\left(\begin{matrix}  \widetilde{\bm \Xi}^{(11)} &  \widetilde{\bm \Xi}^{(12)} \\
            \widetilde{\bm \Xi}^{(21)} &  \widetilde{\bm \Xi}^{(22)}\end{matrix}\right)
\end{align*}
exists.
Also if $i=j=1+\nu_2+\ldots+\nu_l$ and $Re(\lambda_l)=1/2$,
$\nu_l=\nu$, then
\begin{equation}\label{eq7.3}
\begin{array}{rl}
( \widetilde{\bm \Xi}^{(11)})_{ij} & = \frac 1{((\nu-1)!)^2} \frac
1{2\nu-1}
    \big(|\lambda_l|^2\overline{\bm t_{i1}}\bm \Sigma_1\bm t_{j1}^{\prime}
     +\overline{\bm t_{i1}}\bm \Sigma_2\bm t_{j1}^{\prime} \big), \\
(\widetilde{\bm \Xi}^{(22)})_{ij} & = \frac 1{((\nu-1)!)^2} \frac
1{2\nu-1}
    \big( \overline{\bm t_{i1}}\bm \Sigma_1\bm t_{j1}^{\prime}
     +|\lambda_l|^{-2}\overline{\bm t_{i1}}\bm \Sigma_2\bm t_{j1}^{\prime}\big), \\
(\widetilde{\bm \Xi}^{(12)})_{ij} & = \overline{ (\widetilde{\bm
\Xi}^{(21)})_{ji} } =\frac 1{((\nu-1)!)^2} \frac 1{2\nu-1}
    \big( \overline{\lambda}_l\overline{\bm t_{i1}}\bm \Sigma_1\bm t_{j1}^{\prime}
     +\lambda_l^{-1}\overline{\bm t_{i1}}\bm \Sigma_2\bm t_{j1}^{\prime}\big),
\end{array}\end{equation}
and $(\widetilde{\bm \Xi}^{(uv)})_{ij}=0$ for other cases,
$u,v=1,2$. The proof is completed.

\begin{proposition} \label{proposition2.6}Suppose Assumption \ref{assumption1.1} is satisfied. If $\rho<1/2$, then for $i=1,2$,
\begin{equation} \label{eqproposition2.6.1}
\bm G_{ti}=O(t\log\log t)^{1/2})  \;\text{ and }\; \int_0^t\frac{\bm G_{xi}}{x}dx  =O(t\log\log t)^{1/2}) \;a.s. \; t\to \infty.
\end{equation}
If $\rho=1/2$, then for $i=1,2$,
\begin{equation}\label{eqproposition2.6.2}
\begin{array}{rl}
\bm G_{ti}=&O\big((t\log\log\log
t)^{1/2}\log^{v-1/2} t\big)  \; a.s.\;\; t\to \infty, \\
\int_0^t\frac{\bm G_{xi}}{x}dx  =&O\big((t\log\log\log
t)^{1/2}\log^{v-1/2} t\big) \;a.s. \;\; t\to \infty.
\end{array}\end{equation}
\end{proposition}

\noindent {\bf Proof} Let $\bm W_t$ be a $d$-dimensional Brownian
motion, $\bm S_t$  a solution of (Equ1). If $\rho<1/2$, then by
(\ref{eq3.1}) and Lemma \ref{lemma2.7},
\begin{align*}
\|\bm S_t\|=&O((t\log\log t)^{1/2})+\int_0^t \frac{O(x\log\log x)^{1/2}}{x}\big(\frac{t}{x}\big)^{\rho}\log^{v-1}\big(\frac{t}{x}\big)
\\
=&O((t\log\log t)^{1/2})\;\; a.s.
\end{align*}
which implies (\ref{eqproposition2.6.1}).

When $\rho=1/2$, let $\widehat{\bm S}_t$ be  a solution of (Equ1).
Then $\widehat{\bm U}_t=\widehat{\bm S}_t\bm T$ a solution of
(\ref{eq3.20}). It is easily seen that (cf. Bai, Hu and Zhang 2002)
$$
 t^{\lambda_i}\int_1^t \frac{d(\bm W_t \bm t_{ij}^{\prime})}{x^{\lambda_i}}
\overset{a.s.}=
\begin{cases} O\big((t\log\log t)^{1/2}\big)\;\; a.s., & \text{ if } Re(\lambda_i)<1/2 \\
O\big((t\log\log\log t)^{1/2}\log^{1/2}t \big)\;\; a.s., & \text{ if } Re(\lambda_i)=1/2
\end{cases}
$$
as $t\to \infty$. From (\ref{eq3.21}) it follows that, if $Re(\lambda_i)<1/2$,
$$\begin{array}{rl}
\widehat U_{t1}^{(i)}=& O\big((t\log\log t)^{1/2}\big) \;\; a.s.\\
\widehat U_{tj}^{(i)}=& O\big((t\log\log t)^{1/2}\big)
+t^{Re(\lambda_i)}\int_1^t \frac{O((x\log\log
x)^{1/2})}{x^{1+Re(\lambda_i)}}
= O\big((t\log\log t)^{1/2}\big)\;\; a.s.\\
 & j=2,\ldots,v_i; i=1,\ldots,s,
\end{array}$$
and if $Re(\lambda_i)=1/2$,
$$\begin{array}{rl}
\widehat U_{t1}^{(i)}=& O\big((t\log\log\log t)^{1/2}\log^{1/2} t \big)\;\; a.s. \\
\widehat U_{tj}^{(i)}=& O\big((t\log\log\log t)^{1/2}\log^{1/2}
t\big)
+t^{1/2}\int_1^t \frac{O((x\log\log\log x)^{1/2}\log^{j-1-1/2} x )}{x^{1/2}} \\
=& O\big((t\log\log\log t)^{1/2}\log^{j-1/2} t\big)\;\; a.s.\\
& j=2,\ldots,v_i; i=1,\ldots,s.
\end{array}$$
It follows that $\widehat{\bm S}_t=\widehat{\bm U}_t \bm
T^{-1}=O\big((t\log\log\log t)^{1/2}\log^{v-1/2} t\big)$ a.s. Also,
\begin{eqnarray*}
\int_1^t\frac{\widehat{\bm S}_x}x  dx
&= & \int_1^t \frac{O\big((x\log\log\log x)^{1/2}\log^{v-1/2} x\big)}x dx \\
&= & O\big((t\log\log\log t)^{1/2}\log^{v-1/2} t\big)\;\; a.s.
\end{eqnarray*}
So (\ref{eqproposition2.6.2})
 is proved.

\section{Applications.}
\setcounter{equation}{0}

In this section, we give several applications of the approximation theorems.
First, by combing Theorem \ref{theorem1.2} with Proposition \ref{proposition2.4} and Theorem \ref{theorem1.3} with Proposition \ref{proposition2.5} respectively,
  we have the following asymptotic normalities for $(\bm Y_n, \bm
N_n)$.

\begin{theorem}\label{theorem7.1} Under Assumptions \ref{assumption1.1} and
\ref{assumption1.3}, and $\rho<1/2$,
$$n^{-1/2}(\bm Y_n- n\bm v, \bm N_n-n\bm v)\overset{\Cal D}\to
N(\bm 0, \bm \Gamma), $$
where $\bm \Gamma$ is defined in Proposition \ref{proposition2.4}.
\end{theorem}

\begin{theorem}\label{theorem7.2} Suppose  Assumptions  \ref{assumption1.1}
 and \ref{assumption1.3} are satisfied. Further, assume
 (\ref{eq1.18}) is satisfied and $\rho=1/2$. Then
$$n^{-1/2}(\log n)^{1/2-\nu}(\bm Y_n-n\bm v,\bm N_n-n\bm v)\overset{\Cal D}\to
N(0, \widetilde{\bm \Gamma}), $$
where $\widetilde{\bm \Gamma}$ is defined in Proposition \ref{proposition2.5}.
\end{theorem}

Also,  by combining Proposition \ref{proposition2.6}  with  Theorem
\ref{theorem1.2} and Theorem \ref{theorem1.3} respectively,
  we have  we have the following laws of the iterated logarithm.
\begin{theorem}\label{theorem7.3}
Suppose Assumptions \ref{assumption1.1} and \ref{assumption1.3} are
satisfied. If $\rho<1/2$, then
$$\bm Y_n-n\bm v= O\big((n\log\log n)^{1/2}\big)\;\; a.s., $$
$$ \bm N_n- n\bm v= O\big((n\log\log n)^{1/2}\big) \;\; a.s. $$
If (\ref{eq1.18}) is satisfied and $\rho=1/2$, then
$$\bm Y_n-n\bm v = O\big((n\log\log\log n)^{1/2}\log^{\nu-1/2}n \big) \;\; a.s.,$$
$$ \bm N_n-n\bm v = O\big((n\log\log\log n)^{1/2}\log^{\nu-1/2}n \big)\;\; a.s. $$
\end{theorem}

\bigskip
Next, we consider a two-treatment case in which the addition  rule
matrices are denoted by
$$ \bm D_m=\left( \begin{matrix} d_1(\xi_{m1}), & 1- d_1(\xi_{m1}) \\
          1-d_2(\xi_{m2}), & d_2(\xi_{m2}) \end{matrix}\right),$$
where    $(\xi_{11}, \xi_{12}),\ldots,(\xi_{n1},\xi_{n2})$ are
assumed to be i.i.d. random variables with $0\le d_k(\xi_{mk})\le 1$
for $k=1,2$. This is a generalized randomized play-the-winner (RPW)
rule (Bai and Hu, 1999). When $\xi_{m1}$ and $\xi_{m2}$ are
dichotomous and $d_k(x)=x$, then generalized RPW model is the
well-known RPW model proposed by Wei and Durham (1978). In using the
generalized RPW rule, at the stage $m$, if the patient $m$ is
allocated to treatment $1$ and the response $\xi_{m1}$ is observed,
then $d_1(\xi_{m1})$  balls of type $1$ and $1- d_1(\xi_{m1})$ balls
of type $2$ are added to the urn. And, if the patient $m$ is
allocated to treatment $2$ and the response $\xi_{m2}$ is observed,
then $d_2(\xi_{m2})$  balls of type $2$ and $1- d_2(\xi_{m2})$ balls
of type $1$ are added to the urn. It is obvious that  the generating
matrix is
$$ \bm H_m=\bm H=\ep[\bm D_m|\mathscr{F}_{m-1}]=\ep[\bm D_m]=\left(\begin{matrix} p_1, & q_1 \\
          q_2, & p_2 \end{matrix}\right),$$
where $p_k=\ep[d_k(\xi_{mk})]$ and  $q_k=1-p_k$ for $k=1,2$. It is easily checked that Assumptions
\ref{assumption1.1} and \ref{assumption1.3} are  satisfied, and $v_1=q_2/(q_1+q_2)$, $v_2=q_1/(q_1+q_2)$, $\lambda_1=1$, $\rho=\lambda_2=p_1-q_2$.
Denote $\sigma_1^2=v_1v_2=\frac{q_1q_2}{(q_1+q_2)^2 }$ and $\sigma_2^2=\frac{a_1q_2+a_2q_1}{q_1+q_2}$, where $a_k=\Var(d_k(\xi_1))$ for $k=1,2$.
Then
$$\bm \Sigma_1=\sigma_1^2\left(\begin{matrix} 1, & -1 \\ -1, &1
\end{matrix}\right)=\sigma_1^2(1,-1)^{\prime}(1,-1),
\quad \bm \Sigma_2=\sigma_2^2\left(\begin{matrix}
1, & -1 \\ -1, &1 \end{matrix}\right)\sigma_2^2(1,-1)^{\prime}(1,-1),
$$
$$\widetilde{\bm H}=\rho (v_2,-v_1)^{\prime}(1,-1). $$
Further, it is trivial that, if $\bm W_t$ is a Brownian motion with
a variance-covariance matrix $\sigma^2 (1,-1)^{\prime}(1,-1)$, then
$\bm W_t=\sigma(B_t,-B_t)$ where $B_t$ is a standard Brownian
motion. When $\rho<1/2$, multiplying $\bm 1^{\prime}$ in both side
of the equation (Equ1) yields $\bm S_t\bm 1^{\prime}=0$, which
implies $\bm S_t=(S_t,-S_t)$ and $S_t$ is a solution of
$$ S_t=\sigma B_t+\rho \int_0^t\frac{S_x}{x}dx, \;\; S_0=0. $$
It is easily check that
$$ S_t=\sigma t^{\rho}\int_0^t x^{-\rho} d B_x=\sigma B_t+\sigma \rho t^{\rho}\int_0^t B_x x^{-\rho-1} dx \;\; $$
and
$$ \int_0^t \frac{S_x}{x} dx = \sigma   t^{\rho}\int_0^t B_t x^{-\rho-1} dx. $$
Also,
$$
\Var(S_t)= \sigma^2 t^{2\rho}\int_0^t
x^{-2\rho}dx=\frac{\sigma^2}{1-2\rho} t, $$
$$\Var\left\{\int_0^t\frac{S_x}xdx\right\}= \sigma^2 t^{2\rho}
\int_0^t\int_0^t (x\wedge y) x^{-\rho-1}y^{-\rho-1} dx dy
=\frac{2\sigma^2}{(1-2\rho)(1-\rho)}t,  $$
 $$\Cov\{S_t,S_s\}=
\sigma^2\big(\frac{t}{s}\big)^{\rho}\Var(S_s)=\frac{\sigma^2}{1-2\rho}\big(\frac{t}{s}\big)^{\rho}s,
t\ge s, $$
 $$ \Cov\left\{ S_t,\int_0^t\frac{S_x}xdx \right\}=
\int_0^t
\frac{\Cov\{S_t,S_x\}}{x}dx=\frac{\sigma^2}{(1-2\rho)(1-\rho)}t.
$$

 Hence by applying Theorem \ref{theorem1.2} we conclude the
following theorem.

\begin{theorem}For the generalized RPW rule, if $\rho=p_1-q_2<1/2$, then there are two independent standard Brownian motion $B_{t1}$ and $B_{t2}$ such
that for some $\gamma>0$,
\begin{align*}
Y_{n1}-nv_1=&n^{\rho}\int_0^n x^{-\rho} d\big(\rho\sigma_1 B_{x1}+\sigma_2 B_{x2}\big)+o(n^{1/2-\gamma})\;\; a.s.,\\
N_{n1}-n v_1 =& \sigma_1n^{\rho}\int_0^n x^{-\rho} d B_{x1}+\sigma_2
n^{\rho}\int_0^n B_{x2} x^{-\rho-1} dx+o(n^{1/2-\gamma})\;\; a.s.
\end{align*}
and
$$\rho(N_{n1}-n v_1) =n^{\rho}\int_0^n x^{-\rho} d\big(\rho\sigma_1 B_{x1}+\sigma_2 B_{x2}\big)-\sigma_2B_{n2}+o(n^{1/2-\gamma})\;\; a.s. $$
In particular,
$$ n^{1/2}\left(\frac{Y_{n1}}{n}-\frac{q_2}{q_1+q_2},\frac{N_{n1}}{n}-\frac{q_2}{q_1+q_2}\right)\overset{\mathscr{D}}\to N(\bm 0,\bm\Sigma), $$
where $\bm\Sigma=(\sigma_{ij})_{i,j=1}^4$ and
\begin{align*}
\sigma_{11} =& \frac{(p_1-q_2)^2 q_1 q_2+(q_1+q_2)(a_1q_2+a_2 q_1) }{(1-2(p_1-q_2))(q_1+q_2)^2},\\
\sigma_{22}=& \frac{ q_1 q_2+2(a_1q_2+a_2 q_1) }{(1-2(p_1-q_2))(q_1+q_2)^2},\\
\sigma_{12} =\sigma_{21}
                            =&\frac{(p_1-q_2) q_1
q_2+ (a_1q_2+a_2 q_1) }{(1-2(p_1-q_2))(q_1+q_2)^2}.
\end{align*}
\end{theorem}

When $\rho=p_1-q_2=1/2$, by considering the equation (Equ2) and applying
Theorem \ref{theorem1.3} instead, we can define  two independent
standard Brownian motion $B_{t1}$ and $B_{t2}$ such that
\begin{align*}
Y_{n1}-nv_1=&n^{1/2}\int_1^n x^{-1/2} d\big(\frac{1}{2}\sigma_1 B_{x1}+\sigma_2 B_{x2}\big)+O(\sqrt{n})\;\; a.s.,\\
N_{n1}-n v_1 =& \sigma_1n^{1/2}\int_1^n x^{-1/2} d B_{x1}+\sigma_2   n^{1/2}\int_1^n B_{x2} x^{-3/2} dx+O(\sqrt{n})\;\; a.s.
\end{align*}
and
$$\frac{1}{2}(N_{n1}-n v_1) =n^{1/2}\int_1^n x^{-1/2} d\big(\frac{1}{2}\sigma_1 B_{x1}+\sigma_2 B_{x2}\big)-(B_{n2}-n^{1/2} B_{12})+O(\sqrt{n})\;\; a.s. $$
If we denote
$$\widetilde{\sigma}^2=\frac{1}{4}\sigma_1^2+\sigma_2^2=\frac{q_1q_2}{4(q_1+q_2)^2 }+\frac{a_1q_2+a_2q_1}{q_1+q_2}=q_1q_2+2(a_1q_2+a_2q_1).$$
and
$$B(t)=\frac{1}{\widetilde{\sigma}}\int_1^{e^t} x^{-1/2} d\big(\frac{1}{2}\sigma_1 B_{x1}+\sigma_2 B_{x2}\big), $$
it is easily to check that $B(t)$ is a standard Brownian motion. Hence we obtain the following theorem for the case of $\rho=1/2$.

\begin{theorem}For the generalized RPW rule, if $\rho=p_1-q_2=1/2$, then there a  standard Brownian motion $B(t)$ such that
\begin{align*}
Y_{n1}-n\frac{q_2}{q_1+q_2}=&\widetilde{\sigma}n^{1/2}B(\log n)+O(\sqrt{n})\;\; a.s.,\\
N_{n1}-n \frac{q_2}{q_1+q_2} =& 2\widetilde{\sigma}n^{1/2}B(\log n)+\begin{cases} O(\sqrt{n})\;\; & \text{in probability}, \\
 O(n\log\log n)^{1/2} & a.s. \end{cases}
\end{align*}
where $\widetilde{\sigma}^2=q_1q_2+2(q_1q_2+a_2q_1)$. In particular,
$$\limsup_{n\to \infty} \frac{Y_{n1}-nq_2/(q_1+q_2)}{\sqrt{2n(\log  n)(\log\log\log n)}}=\widetilde{\sigma} \;\; a.s., $$
$$\limsup_{n\to \infty} \frac{N_{n1}-nq_2/(q_1+q_2)}{\sqrt{2n(\log  n)(\log\log\log n)}}=2\widetilde{\sigma} \;\; a.s. $$
and
$$ n^{1/2}\left(\frac{Y_{n1}}{n}-\frac{q_2}{q_1+q_2},\frac{N_{n1}}{n}-\frac{q_2}{q_1+q_2}\right)\overset{\mathscr{D}}\to
N\left(\bm 0,\widetilde{\bm\Sigma}\right), $$ where
$\widetilde{\bm\Sigma}=(\widetilde{\sigma}_{ij})_{i,j=1}^4$ and
$\widetilde{\sigma}_{11}=\widetilde{\sigma}^2$,
$\widetilde{\sigma}_{12}=\widetilde{\sigma}_{21}=2\widetilde{\sigma}^2$,
$\widetilde{\sigma}_{22}=4\widetilde{\sigma}^2$.
\end{theorem}


\section{ Proof of the approximation theorems.}
\setcounter{equation}{0}

Define
\begin{equation}\label{eq4.1}\begin{array}{rl}
 \bm M_{n1}  =& \sum_{k=1}^n\{\bm X_k-\ep[\bm X_k|\mathscr{F}_{k-1}]\}
          =: \sum_{k=1}^n \Delta \bm M_{k1},\\
           \bm M_{n2} =&\sum_{m=1}^n\bm X_m(\bm D_m -\ep[\bm D_m |\mathscr{F}_{m-1}])=: \sum_{m=1}^n \Delta \bm M_{m2}.
\end{array} \end{equation}
Recall that $a_n=\bm Y_n \bm 1^{\prime}$ is the total number of
balls in the urn after stage $n$. By (\ref{eq1.1}) we have
\begin{eqnarray}\label{eq4.2}
& & \bm Y_n= \bm Y_0+\sum_{k=1}^n\bm X_k\bm D_k\nonumber
\\
&=&\bm Y_0+\sum_{m=1}^n\big\{\bm X_m(\bm D_m -\ep[\bm D_m
|\mathscr{F}_{m-1}])\nonumber\\
& & \qquad +\big(\bm X_m-\ep[\bm X_m|\mathscr{F}_{m-1}]+\frac {\bm
Y_{m-1}}{a_{m-1}}\big)\bm H+\bm X_m(\bm H_m-\bm H)\big\}
 \nonumber\\
&=& \bm Y_0+\bm M_{n2}+\bm M_{n1}\bm H+\sum_{m=0}^{n-1}\frac {\bm
Y_m}{a_m}\bm H +\sum_{m=1}^n\bm X_m(\bm H_m-\bm H)\nonumber \\
&=&n\bm v+ \bm Y_0+\bm M_{n2}+\bm M_{n1}\bm
H+\sum_{m=0}^{n-1}\big(\frac {\bm
Y_m}{a_m}-\bm v\big)\widetilde{\bm H} +\sum_{m=1}^n\bm X_m(\bm H_m-\bm H)\nonumber \\
&& (\text{ since} \bm Y_m\bm 1^{\prime}=a_m,\;\; \widetilde{\bm H}=\bm H-\bm 1^{\prime} \bm v, \;\; \bm v\widetilde{\bm H}=\bm 0)\nonumber\\
&=&n\bm v+ \bm Y_0+\bm M_{n2}+\bm M_{n1}\bm
H+\sum_{m=1}^{n-1}\big(\frac {\bm
Y_m}{m}-\bm v\big)\widetilde{\bm H} \nonumber\\
& & +\big(\frac {\bm Y_0}{a_0}-\bm v\big)\widetilde{\bm
H}+\sum_{m=1}^{n-1}\frac{m-a_m}{m}\big(\frac {\bm
Y_m}{a_m}-\bm v\big)\widetilde{\bm H}+\sum_{m=1}^n\bm X_m(\bm H_m-\bm H)\nonumber \\
&=:&n\bm v+\bm M_{n2}+\bm M_{n1}\bm H+\sum_{m=1}^{n-1} \frac {\bm
Y_m-m\bm v}{m} \widetilde{\bm H}+\bm R_{n1}+ \bm Y_0,
\end{eqnarray}
where
\begin{eqnarray}\label{eq4.5}
\bm R_{n1} =\big(\frac {\bm Y_0}{a_0}-\bm v\big)\widetilde{\bm
H}+\sum_{m=1}^{n-1}\frac{m-a_m}{m}\big(\frac {\bm Y_m}{a_m}-\bm
v\big)\widetilde{\bm H}+\sum_{m=1}^n\bm X_m(\bm H_m-\bm H).
\end{eqnarray}
Also by (\ref{eq1.2}),
\begin{eqnarray}\label{eq4.3}
\bm N_n&=&\sum_{m=1}^n(\bm X_m-\ep[\bm X_m |\mathscr{F}_{m-1}])
   +\sum_{m=1}^n\ep[\bm X_m |\mathscr{F}_{m-1}]
= \bm M_{n1}+\sum_{m=0}^{n-1}\frac {\bm Y_m}{a_m} \nonumber\\
&=& n\bm v +\bm M_{n1}+\sum_{m=0}^{n-1}\big(\frac{\bm Y_m}{a_m}-\bm v\big)\big(\bm I-\bm 1^{\prime} \bm v\big)\nonumber \\
&=& n\bm v +\bm M_{n1}+\sum_{m=1}^{n-1}\frac{\bm Y_m-m\bm v}{m}
\big(\bm
I-\bm 1^{\prime} \bm v\big)\nonumber \\
& & +\big(\frac{\bm Y_0}{a_0}-\bm v\big)+\sum_{m=1}^{n-1}\frac{ m-a_m}{m}\big(\frac{\bm Y_m}{a_m}-\bm v\big)\big(\bm I-\bm 1^{\prime} \bm v\big)\nonumber \\
  &=&n\bm v +\bm M_{n1}+\sum_{m=1}^{n-1}\frac{\bm Y_m-m\bm v}{m} \big(\bm
I-\bm 1^{\prime} \bm v\big)+\bm R_{n2},
\end{eqnarray}
where
\begin{eqnarray}\label{eq4.4}\bm R_{n2}=\big(\frac{\bm Y_0}{a_0}-\bm v\big)
+\sum_{m=1}^{n-1}\frac{m-a_m}{m}\big(\frac{\bm Y_m}{a_m}-\bm
v\big)\big(\bm I-\bm 1^{\prime} \bm v\big).
\end{eqnarray}

The expansions  given in (\ref{eq4.2}) and (\ref{eq4.3})  are  the
key component in asymptotic analysis of $\bm Y_n$ and $\bm N_n$.
Actually, if we neglect the remainder $\bm R_{n1}$ and replace $\bm
M_{n1}\widetilde{\bm H}+\bm M_{n2}$ by  a Brownian motion $\bm W_n$,
then
$$ \bm Y_n-n\bm v\approx \bm W_n +\sum_{m=1}^{n-1} \frac{\bm
Y_m-m\bm v}{m}\widetilde{\bm H},
$$
which is very similar to the equations (Equ1) or (Equ2).  We will
show $ (\bm Y_n-n\bm v, \bm N_n-n\bm v)$ can be approximated by a
$2d$-dimensional Gaussian process by approximating the martingale
$(\bm M_{n1}, \bm M_{n2})$  to a $2d$-dimensional Brownian motion.
First   show that the remainders $\bm R_{n1}$ and $\bm R_{n2}$ can
be neglected.

\begin{proposition}\label{prop5.1}
Under Assumptions (\ref{assumption1.1}) and (\ref{assumption1.3}),
we have for any $\delta>0$,
\begin{align*}
\bm R_{n1}=&o(n^{\delta})+\sum_{m=1}^n \bm X_m (\bm H_m-\bm
H)=o(n^{1/2-\tau})\;\; a.s., \\
\bm R_{n1}=&o(n^{\delta})\;\; a.s.
\end{align*}
\end{proposition}

To proving this proposition, we need two lemmas, the first one can
be found in Hu and Zhang (2004).

\begin{lemma} \label{lemma5.0} ({\rm  Hu and Zhang (2004)}) If $\bm \Delta \bm Q_n=\Delta \bm
P_n+\bm Q_{n-1} \widetilde{\bm H}/(n-1)$, $n\ge 2$, then
$$ \|\bm Q_n\|=O(\|\bm P_n\|)+\sum_{m=1}^n \frac{O(\|\bm
P_m\|)}{m}\big(n/m)^{\rho}\log^{\nu-1}(n/m). $$
\end{lemma}

\begin{lemma}\label{lemma5.1} Suppose $\sup_m\ep\|\bm D_m\|^2<\infty$. Under Assumptions \ref{assumption1.1},
$$\bm M_n\overset{L_2}=O(n^{1/2})\; \text{ and } \;
\bm M_{n1}\overset{L_2}= O(n^{1/2}), $$
$$\bm M_n\overset{a.s.}=O(n^{1/2+\delta})\; \forall \delta>0
\; \text{ and } \;\bm M_{n1}\overset{a.s.}=O((n\log\log n)^{1/2}),$$
$$
a_n-n\overset{a.s.} =O(n^{1/2+\delta})\; \forall \delta>0. $$
Furthermore, under Assumption \ref{assumption1.3},
$$\bm M_n\overset{a.s.}=O\big((n\log\log n)^{1/2}\big) \; \text{ and }
a_n-n\overset{a.s.}=O\big((n\log\log n)^{1/2}\big). $$
\end{lemma}

\noindent{\bf Proof} Note that $\|\Delta \bm M_{n1}\| \le \|\bm
X_n\|+\ep[\|\bm X_n\| |\mathscr{F}_{n-1}]\le 2$, $ \|\Delta \bm
M_n\|\le \|\bm D_n\|+\ep[\|\bm D_n\| \big |\mathscr{F}_{n-1}]$ and
$a_n=n+\bm Y_0 \bm 1^{\prime}+\bm M_n \bm 1^{\prime}+\sum_{m=1}^n\bm
X_m(\bm H_m-\bm H)\bm 1^{\prime}$. By the properties of martingale, the
results follow easily.

\begin{lemma}\label{lemma5.2} Suppose $\rho\le 1/2$ and    $\sup_m\ep\|\bm D_m\|^2<\infty$.
  Under Assumptions \ref{assumption1.1},
\begin{equation}\label{eq5.1}
\frac {\bm Y_n}{a_n}-\bm v =o(n^{-1/2+\delta}) \; a.s. \text{ for
any } \delta>0. \end{equation}
\end{lemma}

\noindent{\bf Proof} By (\ref{eq4.5})  and Lemma \ref{lemma5.1}, it is obvious  that
$$\|\bm R_{n1} \|\le C
\sum_{m=1}^{n-1}\frac{|m-a_m|}{m}+\sum_{m=1}^n \|\bm H_m-\bm
H\|=o(n^{1/2+\delta})\; a.s. \; \text{ for any  } \delta>0. $$
 From
(\ref{eq4.2}) and Lemma \ref{lemma5.1}, it follows that
$$\bm Y_n- n\bm v
=\sum_{m=1}^{n-1}\frac{\bm Y_m-m \bm v }{m}\widetilde{\bm H}
+o(n^{1/2+\delta})\; a.s. \; $$
 By Lemma \ref{lemma5.0}, it follows that
$$\bm Y_n-n\bm v =
o(n^{1/2+\delta})+\sum_{m=1}^n \frac{o(m^{1/2+\delta})}{m}\big(n/m)^{\rho}\log^{\nu-1}(n/m)=o(n^{1/2+\delta})\; a.s. \;  $$
Hence
$$ \frac {\bm Y_n}{a_n}-\bm v=\frac {\bm Y_n-n\bm v}{n} +\frac{(\bm Y_n-n\bm v)\bm 1^{\prime}}{n}\frac{\bm Y_n}{a_n}=o\big(\frac{n^{1/2+\delta}}{n}\big)
=o(n^{-1/2+\delta})\; a.s. $$
(\ref{eq5.1}) is proved.

Now, we tend to

{\bf Proof of Proposition \ref{prop5.1}.} Notice $\frac{m-a_m}{m} \big(\frac{\bm Y_m}{a_m}-\bm v\big)=o(n^{-1+2\delta})$ a.s. by Lemma \ref{lemma5.2}. The proof is completed by noticing
(\ref{eq4.5}) and (\ref{eq4.4}).

The next result is about the conditional variance-covariance matrix of the $2d$-dimensional   martingale $(\bm M_{n1}, \bm M_{n2})$.
\begin{proposition}\label{proposition3.2} We have
\begin{equation}\label{eqproposition3.2.1}\ep[(\Delta \bm M_{m1})^{\prime}\Delta \bm M_{m2} |\mathscr{F}_{m-1}]=\bm 0
\end{equation}
and under Assumptions \ref{assumption1.1} and
\ref{assumption1.3},
\begin{equation}\label{eqproposition3.2.3} \sum_{m=1}^n  \ep[(\Delta \bm M_{mi})^{\prime} \Delta \bm M_{mi}|\mathscr{F}_{m-1}]
= n  \bm\Sigma_i +o(n^{1-\epsilon})\; a.s. \;\; i=1,2.
\end{equation}
\end{proposition}

\noindent{\bf Proof} (\ref{eqproposition3.2.1}) is trivial. For (\ref{eqproposition3.2.3}), we have
\begin{eqnarray*} &&\ep [\Delta \bm M_{n2}^{\prime}  \Delta \bm M_{n2}|\mathscr{F}_{n-1}] =
 \ep[ (\bm D_n-\bm H_n)^{\prime} diag(\bm X_n)(\bm D_n-\bm H_n)
|\mathscr{F}_{n-1}]\\
&=&
 \ep[ (\bm D_n-\bm H_n)^{\prime} diag(\frac{\bm Y_{n-1} }{a_{n-1}})(\bm D_n-\bm H_n)
|\mathscr{F}_{n-1}]  \\
&=& \ep[ (\bm D_n-\bm H_n)^{\prime} diag(\bm v)(\bm D_n -\bm
H_n)|\mathscr{F}_{n-1}] \\
&&+\ep[ (\bm D_n-\bm H_n)^{\prime} \big(diag(\frac{\bm Y_{n-1}
}{a_{n-1}})-diag(\bm v)\big)(\bm D_n -\bm H_n)|\mathscr{F}_{n-1}] \\
&=& \sum_{q=1}^d v_q\bm V_{nq}
+\sum_{q=1}^d (\frac{Y_{n-1,q}}{a_{n-1}}- v_q)\bm V_{nq}.
\end{eqnarray*}
 Under Assumptions \ref{assumption1.1},
\ref{assumption1.3}, by  Lemma \ref{lemma5.2} we have
\begin{eqnarray*}
\sum_{m=1}^n \ep [\Delta \bm M_{m2}^{\prime}  \Delta \bm M_{m2}|\mathscr{F}_{m-1}]= n \bm\Sigma_2+o(n^{1-\epsilon})\;\; a.s.
\end{eqnarray*}
Also
\begin{eqnarray*}
&&\ep [\Delta \bm M_{n1}^{\prime}
 \Delta \bm M_{n1}|\mathscr{F}_{n-1}]\\
&=&  \ep [ \bm X_n^{\prime} \bm X_n |\mathscr{F}_{n-1}]
-\big(\ep[\bm X_n|\mathscr{F}_{n-1}]\big)^{\prime}\ep[\bm X_n|\mathscr{F}_{n-1}] \\
&=& \ep[  diag(\bm X_n) |\mathscr{F}_{n-1}]
 -\frac{ Y_{n-1}^{\prime}}{a_{n-1}}\frac{ Y_{n-1} }{a_{n-1}}
= diag(\frac{\bm Y_{n-1} }{a_{n-1}})
 -\frac{ Y_{n-1}^{\prime}}{a_{n-1}}\frac{ Y_{n-1} }{a_{n-1}} \\
&=&diag(\bm v)-\bm v^{\prime}\bm v+ o(n^{-1/2+\delta})\;\; a.s.
\end{eqnarray*}
(\ref{eqproposition3.2.3}) is proved.


\bigskip

{\bf Proof of Theorem \ref{theorem1.2}.}  Suppose that Assumption \ref{assumption1.3}
is satisfied. According to (\ref{eqproposition3.2.3}),
\begin{eqnarray*}
&&\sum_{m=1}^n  \ep[(\Delta \bm M_{m1},\Delta \bm
M_{m2})^{\prime}(\Delta \bm M_{m1},\Delta \bm
M_{m2})|\mathscr{F}_{m-1}]=n\; diag(\bm
\Sigma_1,\bm\Sigma_2)+o(n^{1-\epsilon})\;\; a.s.
\end{eqnarray*}
It follows from Theorem 1.3 of Zhang (2004)  that, there exist two
standard $d$-dimensional Brownian motions $\bm B_{t1}$ and $\bm
B_{t2}$ for which
\begin{equation}\label{eq6.17}
 (\bm M_{n1},
\bm M_{n2})-(\bm B_{n1}\bm\Sigma_1^{1/2}, \bm B_{n2}\bm\Sigma_2^{1/2})=o(n^{1/2-\gamma})\quad a.s.
\end{equation} Here $\gamma>0$ depends only on $d$ and $\epsilon$. Without
loss of generality, we assume $\gamma<\epsilon/3$.

Now let $\bm G_{ti}=\text{Solut}(\text{Equ1}, \bm
B_{ti}\Sigma_i^{1/2})$ ($i=1,2$). Then by Proposition
\ref{proposition2.3},
$$ \int_0^n \frac{\bm G_{xi}}{x}dx=\sum_{m=1}^{n-1}\frac{\bm G_{mi}}{m}+O(1) \;a.s. \;\; i=1,2. $$
Write $\bm G_t=\bm G_{t2}\bm\Sigma_2^{1/2}+\bm G_{t1}\bm \Sigma_1^{1/2}\bm H$.
Combining the above equality with (\ref{eq4.2}), (\ref{eq6.17}) and Proposition \ref{prop5.1} yields
$$\bm Y_n-n\bm v-\bm G_n=\sum_{m=1}^{n-1}
\frac {\bm Y_m-m\bm v-\bm G_m}{m}\widetilde{\bm H}
+o(n^{1/2-\tau\wedge\gamma})\;\; a.s. $$ By Proposition
\ref{lemma5.0},
$$ \bm Y_n-n\bm v-\bm G_n=o(n^{1/2-\tau\wedge\gamma})+\sum_{m=1}^n\frac{o(m^{1/2-\tau\wedge\gamma})}{m}(n/m)^{\rho}\log^{v-1}(n/m)=o(n^{1/2-\tau\wedge\gamma})\;\; a.s. $$
Finally, combining the above equality with (\ref{eq4.3}),
(\ref{eq6.17}) and Proposition \ref{prop5.1} yields
\begin{eqnarray*}
\bm N_n- n\bm v&=& \bm M_{n1}
  +\sum_{m=0}^{n-1}\frac {\bm Y_m- \ep\bm Y_m}{m}(\bm I-\bm 1^{\prime}\bm v)
+o(n^{\delta}) \\
&=&  \bm B_{n1}\bm \Sigma_1^{1/2}
  +\sum_{m=1}^{n-1}\frac {\bm G_m}{m}(\bm I-\bm 1^{\prime}\bm v)
+o(n^{1/2-\tau\wedge\gamma}) \\
&=&  \bm B_{n1}\bm \Sigma_1^{1/2}
  +\int_0^n\frac {\bm G_x}{x}\;dx(\bm I-\bm 1^{\prime}\bm v)
+o(n^{1/2-\tau\wedge\gamma}) \\
&=& \bm G_{n1} +\int_0^n\frac {\bm G_{x2}}{x}\;dx(\bm I-\bm
1^{\prime}\bm v) +o(n^{1/2-\tau\wedge\gamma})\;\; a.s.
\end{eqnarray*}
The proof is now completed.

\bigskip

{\bf Proof of Theorem \ref{theorem1.3}.}  (\ref{eq6.17}) remains
true. Let $\widehat{\bm G}_{ti}=\text{Solut}(\text{Equ2}, \bm
B_{ti}\Sigma_i^{1/2})$ ($i=1,2$). Then by Proposition
\ref{proposition2.3},
$$ \int_1^n \frac{\widehat{\bm G}_{xi}}{x}dx=\sum_{m=1}^{n-1}\frac{\widehat{\bm G}_{mi}}{m}+O(1) \;a.s. \;\; i=1,2. $$
Write $\widehat{\bm G}_t=\widehat{\bm G}_{t2}\bm\Sigma_2^{1/2}+\widehat{\bm G}_{t1}\bm \Sigma_1^{1/2}\bm H$.
Combining the above equality with (\ref{eq4.2}), (\ref{eq6.17}) and Proposition \ref{prop5.1} yields
$$\bm Y_n-n\bm v-\widehat{\bm G}_n=\sum_{m=1}^{n-1}
\frac {\bm Y_m-m\bm v-\widehat{\bm G}_m}{m}\widetilde{\bm H} +o(n^{1/2-\gamma})+\sum_{m=1}^n \bm X_m(\bm X_m-\bm H)\;\; a.s. $$
By Proposition \ref{lemma5.0},
\begin{align*}
\bm Y_n- & n\bm v-\widehat{\bm G}_n=o(n^{1/2-\gamma})
+\sum_{m=1}^n\frac{o(m^{1/2-\gamma})}{m}(n/m)^{1/2}\log^{v-1}(n/m)\\
& +O(\sum_{m=1}^n\|\bm H_m-\bm H\|)
+\sum_{m=1}^n\frac{O(\sum_{j=1}^m\|\bm H_m-\bm H\|)}{m}(n/m)^{1/2}\log^{v-1}(n/m)\\
=& o(n^{1/2})+ O(1)\sum_{j=1}^n\frac{\|\bm H_j-\bm H\|}{j^{1/2}} n^{1/2}\log^{\nu-1}n=O\big(n^{1/2}\log^{\nu-1}n\big)\;\; a.s.
\end{align*}
Finally, combining the above equality with (\ref{eq4.3}), (\ref{eq6.17}) and Proposition \ref{prop5.1} yields
\begin{eqnarray*}
\bm N_n- n\bm v&=& \bm M_{n1}
  +\sum_{m=0}^{n-1}\frac {\bm Y_m- \ep\bm Y_m}{m}(\bm I-\bm 1^{\prime}\bm v)
+o(n^{\delta}) \\
&=&  \bm B_{n1}\bm \Sigma_1^{1/2}
  +\sum_{m=1}^{n-1}\frac {\widehat{\bm G}_m}{m}(\bm I-\bm 1^{\prime}\bm v)
+O\big(n^{1/2}\log^{\nu-1}n\big) \\
&=&  \bm B_{n1}\bm \Sigma_1^{1/2}
  +\int_1^n\frac {\widehat{\bm G}_x}{x}\;dx(\bm I-\bm 1^{\prime}\bm v)
+O\big(n^{1/2}\log^{\nu-1}n\big) \\
&=& \widehat{\bm G}_{n1} +\int_1^n\frac {\widehat{\bm G}_{x2}}{x}\;dx(\bm
I-\bm 1^{\prime}\bm v) +O\big(n^{1/2}\log^{\nu-1}n\big)\;\; a.s.
\end{eqnarray*}
The proof is now completed.

\pagebreak \baselineskip 14pt \linewidth 6.8in \oddsidemargin=-0.1in
\begin{center}
REFERENCES
\end{center}
\baselineskip = 14pt
\begin{verse}
\vspace{-0.28in} \hspace{-0.5in} {\small \baselineskip 14pt
\item
Andersen, J., Faries, D. and Tamura, R. N. (1994).  Randomized
play-the-winner design for multi-arm clinical trials.  {\em
Communications in Statistics, Theory and Methods} {\bf 23} 309-323.
\item
Athreya, K. B. and Karlin, S. (1967). Limit theorems for the split
times of branching processes. {\em Journal of Mathematics and
Mechanics}  {\bf 17} 257-277.
\item
Athreya, K. B. and Karlin, S. (1968). Embedding of urn schemes into
continuous time branching processes and related limit theorems. {\em
Ann. Math. Statist.} {\bf 39} 1801-1817.
\item
Bai, Z. D. and Hu, F. (1999). Asymptotic theorem for urn models with
nonhomogeneous generating matrices. {\em Stochastic Process. Appl.}
{\bf 80} 87-101.
\item
Bai, Z. D. and Hu, F. (2005). Asymptotics in randomized urn models.
{\em Ann. Appl. Probab.} {\bf 15}  914-940.
\item
Bai, Z. D., Hu, F. and Shen, L. (2002). An adaptive design for
multi-arm clinical trials, {\em J. Multi. Anal.}. {\bf 81} 1-18.
\item
 Bai, Z. D., Hu, F. and Zhang, L. X.  (2002).
The Gaussian approximation theorems for urn models and their
applications. {\em Ann. Appl. Probab.} {\bf 12} 1149-1173.
\item
Beggs, A. W. (2005). On the convergence of reinforcement learning.
{\em Journal of Economic Theory} To appear.
\item
Bena\"{i}m M., Schreiber S. J. and Tarr\'{e}s, P. (2004).
Generalized urn models of evolutionary processes. {\em The Annals of
Applied Probability} {\bf 14} 1455-1478.
\item
Eggenberger, F. and P\'{o}lya, G. (1923). \"{U}ber die statistik
verketetter vorg\"{a}nge. {\em Zeitschrift f\"{u}r Angewandte
Mathematik und Mechanik} {\bf 1} 279-289.
\item
Erev, I. and Roth, A. (1998). Predicting how people play games:
reinforcement learning in experimental games with unique, mixed
strategy equilibria.  {\em Amer. Econ. Rev.} {\bf 88} 848-881.
\item Flournoy, N. and Rosenberger, W. F., eds. (1995). {\em
Adaptive Designs}. Hayward, Institute of Mathematical Statistics.
\item   Gouet, R. (1993).
Martingale functional central limit theorems for a generalized
P\'olya urn. {\em Ann. Probab.} {\bf 21} 1624-1639.
\item
Hall, P. and Heyde, C. C. (1980). {\it Martingale Limit Theory and
its Applications}. Academic Press, London.
\item
Hu, F. and Rosenberger, W. F. (2003). Analysis of time trends in
adaptive designs with application to a neurophysiology experiment.
{\em Statist. Med.}  {\bf 19} 2067-2075.
\item   Hu, F.  and   Rosenberger, W. F.  (2006).
{\em The Theory of Response-Adaptive Randomization in Clinical
Trials}. John Wiley and Sons, Inc., New York.
\item
Hu, F. and Zhang, L. X. (2004). Asymptotic properties of doubly
adaptive biased coin designs for multi-treatment clinical trials.
{\em Ann. Statist.} {\bf 32} 268-301.
\item
Janson, S. (2004). Functional limit theorems for multitype branching
processes and generalized P\'olya urns. {\em Stochastic Process.
Appl.} {\bf 110} 177-245.
\item
Johnson, N. L. and Kotz, S. (1977). {\em Urn Models and Their
Applications}. Wiley, New York.
\item
Kotz, S. and Balakrishnan, N. (1997). Advances in urn models during
the past two decades, In {\em Advances in Combinatorial Methods and
Applications to Probability and Statistics} (Eds., Balakrishnan, N.)
Birkh\"{a}user, Boston.
\item
Martin, C. F. and Ho, Y. C. (2002). Value of information in the
Polya urn process. {\em Information Sciences} {\bf 147} 65-90.
\item
Robbins, H. (1952).
 Some aspects of the sequential design of
experiments. {\em Bull. Amer. Math. Soc.} {\bf 58} 527-535.
\item   Rosenberger, W. F. (1996). New directions in adaptive
designs. {\it Statist. Sci.}  {\bf 11} 137-149.
\item
Rosenberger, W. F. and  Lachin, J. M.  (2002). {\em Randomization in
Clinical Trials:  Theory and Practice.} Wiley,  New York.
\item
Smythe, R. T. (1996). Central limit theorems for urn models. {\em
Stochastic Process. Appl.} {\bf 65} 115-137.
\item
Smythe, R. T. and Rosenberger, W. F. (1995). Play-the-winner
designs, generalized P\'olya urns, and Markov branching processes.
In {\em Adaptive Designs} (Flournoy, N. and Rosenberger, W. F. eds.)
Hayward, CA: Institute of Mathematical Statistics, 13-22.
\item
Thompson, W. R. (1933).
 On the likelihood that one
unknown probability exceeds another in view of the evidence of the
two samples. {\em Biometrika} {\bf 25} 275-294.
\item
Wei, L. J. (1979). The generalized P\'olya's urn design for
sequential medical trials. {\it Ann. Statist.} {\bf 7} 291-296.
\item
Wei, L. J. and Durham, S. (1978). The randomized pay-the-winner rule
in medical trials. {\it J. Amer. Statist. Assoc.} {\bf 73} 840-843.
\item
Windrum, P. (2004). Leveraging technological externalities in
complex technologies: Microsoft's exploitation of standards, in the
browser wars. {\em Research Policy} {\bf 33} 385-394.}
\item
Zhang, L. X. (2004). Strong approximations of martingale vectors and its
applications in Markov-Chain adaptive designs. {\em Acta Math. Appl.
Sinica, English Series} {\bf 20}(2) 337-352
\end{verse}

\bigskip

\begin{center}
{\sc \footnotesize{

\begin{tabular}{lcl}
L-X. Zhang & & F. Hu \\ [-0.30cm] Department of Mathematics & &
Department of Statistics \\ [-0.30cm] Zhejiang University & &
University of Virginia \\ [-0.30cm] Hangzhou 310027 & & Halsey Hall,
Charlottesville \\ [-0.30cm] PR China & & Virginia 22904-4135, USA \\
[-0.30cm] E-Mail: {\rm stazlx@zju.edu.cn} & & E-Mail: {\rm
fh6e@virginia.edu}
\end{tabular}
}}

\end{center}

\end{document}